\documentclass{birkjour}
\usepackage[latin1]{inputenc}
\usepackage{qcftw}

\usepackage{latexsym}
\usepackage{amsmath}
\usepackage{amsfonts}
\usepackage{amssymb}
\usepackage{amsthm}
\usepackage{color} 
\usepackage{verbatim}
\usepackage{subfigure}
 
\newcommand{\bomega}{\mbox{\boldmath $\omega$}}
\newcommand{\bmu}{\mbox{\boldmath $\mu$}}
\newcommand{\bvect}[1]{\mbox{\boldmath $#1$}} 
\newcommand{\sandwich}{(\,)} 
\newcommand{\be}{\begin{equation}}
\newcommand{\ee}{\end{equation}}	

\makeindex

\begin{document}

\title[Orthogonal 2D \change{P}lanes \change{S}plit]%
      {The Orthogonal 2D Planes Split of\\ Quaternions and Steerable Quaternion\\ Fourier Transformations}

\author[E. Hitzer]{Eckhard Hitzer}
\address{%
College of Liberal Arts, Department of Material Science,\\ 
International Christian University,\\
181-8585 Tokyo, Japan}
\email{hitzer@icu.ac.jp}

\author[Sangwine]{Stephen J. Sangwine}

\address{School of Computer Science and Electronic Engineering,\\
         University of Essex, Wivenhoe Park, Colchester, CO4 3SQ, UK.}
         
\email{sjs@essex.ac.uk}

\keywords{quaternion signals, orthogonal 2D planes split, quaternion Fourier transformations, steerable transforms, geometric interpretation, fast implementations}

\begin{abstract}
The two-sided quaternionic Fourier transformation (QFT) was introduced in \cite{Ell:1993}
for the analysis of 2D linear time-invariant partial-differential systems.
In further theoretical investigations \cite{10.1007/s00006-007-0037-8, EH:DirUP_QFT}
a special split of quaternions was introduced,
then called $\pm$split.
In the current \change{chapter} we analyze this split further,
interpret it geometrically as \change{an} \emph{orthogonal 2D planes split} (OPS),
and generalize it to a freely steerable split of $\H$ into two orthogonal 2D analysis planes\index{quaternions!orthogonal 2D planes split}.
The new general form of the OPS split allows us to find new geometric
interpretations for the action of the QFT on the signal.
The second major result of this work is
a variety of \emph{new steerable forms} of the QFT\index{quaternion Fourier transform!steerable}, their geometric interpretation,
and for each form\change{,} OPS split theorems,
which allow fast and efficient numerical implementation with standard FFT software.
\end{abstract}

\maketitle


\section{Introduction}

The two-sided quaternionic Fourier transformation (QFT) was introduced in \cite{Ell:1993}
for the analysis of 2D linear time-invariant partial-differential systems.
Subsequently it has been applied in many fields,
including colour image processing\index{quaternions!colour image!processing}\index{colour image!processing} \cite{10.1049/el:19961331}.
This led to further theoretical investigations \cite{10.1007/s00006-007-0037-8, EH:DirUP_QFT},
where a special split of quaternions was introduced, then called \change{the} $\pm$split.
An interesting physical consequence was that this split resulted \change{in a} left and right travelling multivector wave packet analysis, when generalizing the QFT
to a full spacetime Fourier transform (SFT). In the current \change{chapter} we investigate this split further, interpret it geometrically and generalize it to a
\emph{freely steerable}\footnote{Compare Section \ref{sec:determinesteer},            in particular Theorem \ref{th:abdetfg}.}
split of $\H$ into two orthogonal 2D analysis planes.
For reasons to become obvious we prefer to call it from now on the
\emph{orthogonal 2D planes split} (OPS).

The general form of the OPS split allows us to find new geometric interpretations for the action of the QFT on the signal. The second major result of this work \change{is} a variety of new forms of the QFT, their detailed geometric interpretation, and for each form, \change{an} OPS split \change{theorem}, which allow\change{s} fast and efficient numerical implementation with standard FFT software.
A preliminary formal investigation of these new OPS-QFTs can be found in \cite{EH:OPS-QFT}.

The chapter is organized as follows. We first introduce in Section \ref{sc:pmsplit} several properties of quaternions together with a brief review of the $\pm$-\emph{split} of \cite{10.1007/s00006-007-0037-8, EH:DirUP_QFT}. In Section \ref{sc:OPSfg} we generalize this split to a \emph{freely steerable orthogonal 2D planes split} (OPS) of quaternions $\H$. In Section \ref{sc:NewQFT} we use the general OPS of Section \ref{sc:OPSfg} to generalize the two sided QFT to a \emph{new two sided QFT} with freely \emph{steerable analysis planes}, complete with a detailed \emph{local geometric transformation interpretation}. The geometric interpretation of the OPS in Section \ref{sc:OPSfg} further allows the construction of a new type of \emph{steerable QFT with a direct phase angle interpretation}. In Section \ref{sc:conjQFT} we finally investigate \emph{new steerable QFTs involving quaternion conjugation}. Their local geometric interpretation crucially relies on the notion of \emph{4D rotary reflections}.

\section{Orthogonal Planes Split of Quaternions with Two\\
         Orthonormal Pure Unit Quaternions\label{sc:pmsplit}}

Gauss, Rodrigues and Hamilton's four-dimensional (4D) quaternion algebra\index{quaternions!algebra} $\H$ is defined over $\R$ 
with three imaginary units:

\begin{gather}
 \i \j = -\j \i = \k, \qquad
 \j \k = -\k \j = \i, \qquad
 \k \i = -\i \k = \j,
 \nonumber \\
 \i^2=\j^2=\k^2=\i \j \k = -1.
\label{eq:quat}
\end{gather}
Every quaternion can be written explicitly as
\begin{equation}
  q=q_r + q_i \i + q_j \j + q_k \k \in \H, \quad 
  q_r,q_i, q_j, q_k \in \R,
  \label{eq:aquat}
\end{equation}
and has a \emph{quaternion conjugate}
(equivalent\footnote{\change{This may be important in generalisations of the QFT,
such as to a space-time Fourier transform in \cite{10.1007/s00006-007-0037-8}, or a general two-sided
Clifford Fourier transform in \cite{EH:Agacse2012}.}}
to Clifford conjugation in \change{$\clifford{3,0}^+$ and $\clifford{0,2}$})
\begin{equation}
  \qconjugate{q} = q_r - q_i \i - q_j \j - q_k \k,
  \quad 
  \qconjugate{pq} = \qconjugate{q}\,\qconjugate{p},
  \label{eq:quatconj}
\end{equation} 
which leaves the scalar part $q_r$ unchanged.
This leads to the \emph{norm} of $q\in\H$
\begin{equation}
  \modulus{q} = \sqrt{q\qconjugate{q}} = \sqrt{q_r^2+q_i^2+q_j^2+q_k^2},
  \qquad
  \modulus{p q} = \modulus{p}\modulus{q}.
\end{equation}
The part $\Vector{q} = q - q_r = \frac{1}{2}(q-\qconjugate{q}) = q_i \i + q_j \j + q_k \k$ is called a
\emph{pure} quaternion, and it squares to the negative number $-(q_i^2+q_j^2+q_k^2)$.
Every unit quaternion (\ie $\modulus{q}=1$) can be written as:
\begin{align}
q &= q_r + q_i \i + q_j \j + q_k \k
  = q_r + {\sqrt{q_i^2+q_j^2+q_k^2}}\,{\bmu(q)}
  \nonumber \\
  &= \cos\alpha + {\bmu(q)}\sin\alpha 
  = e^{\alpha\,{\bmu(q)}},\label{eq:unit}\\
  \intertext{where}
  &\cos \alpha = q_r, \qquad 
  \sin \alpha = \sqrt{q_i^2+q_j^2+q_k^2},
  \nonumber\\
  {\bmu(q)}&= \frac{\Vector{q}}{\modulus{q}} = \frac{q_i \i + q_j \j + q_k\k}{\sqrt{q_i^2+q_j^2+q_k^2}},
  \qquad
  \text{and}\qquad {\bmu(q)}^2 = -1.
  \label{eq:unitq}
\end{align}
The \emph{inverse} of a non-zero quaternion is
\begin{equation}
  \inverse{q} = \frac{\qconjugate{q}}{\modulus{q}^2} = \frac{\qconjugate{q}}{q\qconjugate{q}}.
\end{equation} 
The \emph{scalar part} of a quaternion is defined as
\begin{equation}
  \Scalar{q} = q_r = \frac{1}{2}(q+\qconjugate{q}),
\end{equation}
with \emph{symmetries}\index{quaternions!scalar part!symmetries}
\begin{equation}
  \Scalar{pq} = \Scalar{qp} = p_rq_r - p_iq_i - p_jq_j - p_kq_k,
  \quad 
  \Scalar{q} = \Scalar{\qconjugate{q}}, 
  \quad
  \forall p,q \in \H,
  \label{eq:Scsymm}
\end{equation} 
and \emph{linearity}
\begin{equation}
  \Scalar{\alpha p+ \beta q}
  = \alpha \Scalar{p} + \beta \Scalar{q} 
  = \alpha p_r + \beta q_r,
  \quad 
  \forall p,q \in \H, \,\,\, \alpha, \beta \in \R.
  \label{eq:Sclin}
\end{equation} 
The scalar part and the quaternion conjugate allow the definition of
the $\R^4$ \emph{inner product}\footnote{\change{Note that we do not use the notation $p\cdot q$, which is unconventional for full quaternions.}}
of two quaternions $p,q$ as
\begin{equation}
  \Scalar{p\qconjugate{q}}
  = p_rq_r + p_iq_i + p_jq_j + p_kq_k \in \R . 
  \label{eq:4Dinnp}
\end{equation}
\begin{defn}[Orthogonality of quaternions\index{quaternions!orthogonality}]\label{df:qorth}
Two quaternions $p,q \in \H$ are \emph{orthogonal} $p\perp q$, if and only if
\change{the inner product} $\Scalar{p\qconjugate{q}}=0$.
\end{defn}
The \emph{orthogonal\/\footnote{\change{Compare Lemma \ref{lm:smxprod}.}} 2D planes split}
(OPS) of quaternions with respect to the orthonormal pure unit quaternions $\i, \j$ \cite{10.1007/s00006-007-0037-8, EH:DirUP_QFT} is defined by
\begin{gather}
  q = q_+ + q_-, \quad q_{\pm} = \frac{1}{2}(q\pm \i q \j).
  \label{eq:pmform}
\end{gather}
Explicitly in real components
$q_r,q_i, q_j, q_k \in \R$ using \eqref{eq:quat} we get
\begin{equation}
  q_{\pm} = \{q_r\pm q_k + \i(q_i\mp q_j)\}\frac{1\pm \k}{2}
       = \frac{1\pm \k}{2} \{q_r\pm q_k + \j(q_j\mp q_i)\}.
  \label{eq:qpm}
\end{equation}
This leads to the following new Pythagorean \emph{modulus identity\index{quaternions!modulus identity}} \cite{EH:DirUP_QFT}
\begin{lem}[Modulus identity]
  \label{lm:modid}
For $q \in \H$
\begin{equation}
  \label{eq:modid}
  \modulus{q}^2 = \modulus{q_-}^2 + \modulus{q_+}^2.
\end{equation}
\end{lem}

\begin{lem}[Orthogonality of OPS split parts\index{quaternion split!orthogonality}]\label{lm:smxprod}
  Given any two quaternions $p,q \in \H$ and applying the OPS of \eqref{eq:pmform}
  the resulting parts are orthogonal 
\begin{equation}
   \Scalar{p_+\qconjugate{q_-}} = 0, \qquad 
   \Scalar{p_-\qconjugate{q_+}} = 0 ,
\end{equation}  
\ie $p_+ \perp q_-$ and $p_- \perp q_+$.  
\end{lem}
In Lemma \ref{lm:smxprod} (proved in \cite{EH:DirUP_QFT}) the second identity follows from the first by
$\Scalar{\qconjugate{x}} = \Scalar{x}, \, \forall x \in \H$, and  
$\qconjugate{p_-\qconjugate{q_+}} = q_+\qconjugate{p_-}$. 

It is evident, that instead of $\i, \j$\change{,} any pair of orthonormal pure quaternions can be used to produce an analogous split. This is a first indication, that the OPS of \eqref{eq:pmform} is in fact \emph{steerable\index{quaternion split!steerable}}. We observe, that $\i q \j = q_+ - q_-$, \ie under the map $\i\sandwich\j$ the $q_+$ part is invariant, the $q_-$ part changes sign. Both parts are according to \eqref{eq:qpm} two-dimensional, and by Lemma \ref{lm:smxprod} they span two completely orthogonal planes. The $q_+$ plane is spanned by the orthogonal quaternions $\{\i - \j, 1+\i\j=1+\k\}$, whereas the  $q_-$ plane is \eg spanned by $\{\i + \j, 1-\i\j=1-\k\}$, \ie we have the two 2D subspace \change{bases}\index{quaternions!plane subspace bases}
\be 
  q_+\text{-basis: } \{\i - \j, 1+\i\j=1+\k\},\qquad
  q_-\text{-basis: } \{\i + \j, 1-\i\j=1-\k\}.
  \label{eq:q-q+basisij}
\ee  
Note that all basis vectors of \eqref{eq:q-q+basisij}
\be
  \{\i - \j, 1+\i\j,\i + \j, 1-\i\j\}
  \label{eq:ijsplitbasis}
\ee
together form an orthogonal basis\index{quaternions!orthogonal basis} of $\H$ interpreted as $\R^4$. 

The map $\i\sandwich\j$ rotates the $q_-$ plane by $180^{\circ}$\index{quaternions!half-turn rotation} around the 2D $q_+$ axis plane. Note that in agreement with its geometric interpretation, the map $\i\sandwich\j$ is an \emph{involution}\index{quaternions!involution}, because applying it twice leads to identity
\be
  \i(\i q \j)\j = \i^2 q \j^2 = (-1)^2 q = q. 
\ee

\section{General Orthogonal 2D Planes Split}\label{sc:OPSfg}
We will study generalizations of the OPS split\index{quaternions!orthogonal 2D planes split!general} by replacing $\i, \j$ by arbitrary unit quaternions $f, g$. Even with this generalization, the map $f\sandwich g$ continues to be an involution, because $f^2 q g^2 = (-1)^2q = q$. For clarity we study the cases $f\neq \pm g$, and $f=g$ separately, though they have a lot in common, and do not always need to be distinguished in \change{specific} applications.

\subsection{Orthogonal 2D Planes Split using Two Linearly Independent Pure Unit Quaternions}
Our result is now, that all these properties hold, even if in the above considerations the pair $\i, \j$ is replaced by an arbitrary pair of linearly independent nonorthogonal pure quaternions $f,g$, $f^2=g^2=-1, f\neq \pm g$. The OPS is then \emph{re-defined} with respect to the linearly independent pure unit quaternions $f, g$ as 
\begin{equation}
  \label{eq:opsfqg}
  q_{\pm} = \frac{1}{2}(q \pm f q g).
\end{equation} 
Equation \eqref{eq:pmform} is a special case with $f=\i, g=\j$. We  observe from \eqref{eq:opsfqg}, that $f q g = q_+ - q_-$, \ie under the map $f\sandwich g$ the $q_+$ part is invariant, but the $q_-$ part changes sign
\begin{equation}
  fq_{\pm}g 
  = \frac{1}{2}(fqg \pm f^2 q g^2)
  = \frac{1}{2}(fqg \pm q )
  = \pm \frac{1}{2}(q \pm f q g) 
  = \pm q_{\pm}.
  \label{eq:fgrotqm}
\end{equation} 
We now show that even for \eqref{eq:opsfqg} both parts are two-dimensional, and span two completely orthogonal planes. The $q_+$ plane is spanned\footnote{\change{For $f=\i, g=\j$ this is in agreement with \eqref{eq:qpm} and \eqref{eq:q-q+basisij}!}} by the \emph{orthogonal} pair of quaternions $\{f-g, 1+fg\}$:
\begin{gather} 
  \Scalar{(f-g)\qconjugate{(1+fg)}}
  = \Scalar{(f-g)(1+(-g)(-f))}
  \nonumber \\
  = \Scalar{f +fgf -g -g^2f} 
  \stackrel{\eqref{eq:Scsymm}}{=} 
  \Scalar{f +f^2g -g +f} = 2 \Scalar{f-g}=0,
  \label{eq:plusortho}
\end{gather}
whereas the  $q_-$ plane is \eg spanned by $\{f+g, 1-fg\}$. The quaternions $f+g, 1-fg$ can be proved to be mutually \emph{orthogonal} by simply replacing $g\rightarrow -g$ in \eqref{eq:plusortho}. Note that we have
\begin{gather} 
  f(f-g)g = f^2g - fg^2 = -g+f = f-g,
  \nonumber \\
  f(1+fg)g = fg + f^2g^2 = fg+1 = 1+fg,
  \label{eq:fgplus}
\end{gather}
as well as
\begin{gather} 
  f(f+g)g = f^2g + fg^2 = -g-f = -(f+g),
  \nonumber \\
  f(1-fg)g = fg - f^2g^2 = fg-1 = -(1-fg).
  \label{eq:fgminus}
\end{gather}
We now want to generalize Lemma \ref{lm:smxprod}. 
\begin{lem}[Orthogonality of two OPS planes]\label{lm:OPSortho}\index{quaternions!orthogonality of OPS planes}
  Given any two quaternions $q,p \in \H$ and applying the OPS \eqref{eq:opsfqg} with respect to two linearly independent pure unit quaternions $f,g$ we get zero for the
  scalar part of the mixed products
  \begin{equation}
     \Scalar{p_+\qconjugate{q_-}} = 0, 
     \qquad \Scalar{p_-\qconjugate{q_+}} = 0 .
  \end{equation}
\end{lem}
We prove the first identity, the second follows from $\Scalar{x} = \Scalar{\qconjugate{x}}$. 
\begin{gather}
  \Scalar{p_+\qconjugate{q_-}} 
  = \frac{1}{4}\Scalar{(p+fpg)(\qconjugate{q}-g\qconjugate{q}f)}
  = \frac{1}{4}\Scalar{p\qconjugate{q}-fpgg\qconjugate{q}f+fpg\qconjugate{q}-pg\qconjugate{q}f}
  \nonumber \\
  \stackrel{\eqref{eq:Sclin},\eqref{eq:Scsymm}}{=} 
    \frac{1}{4}\Scalar{p\qconjugate{q}-p\qconjugate{q}+pg\qconjugate{q}f - pg\qconjugate{q}f} = 0.
  \label{eq:prooffgortho}
\end{gather}

Thus the set 
\be
  \{f-g, 1+fg, f+g, 1-fg\}
  \label{eq:fgsplitbasis}    
\ee
forms a 4D orthogonal basis\index{quaternions!orthogonal basis} of $\H$ interpreted by \eqref{eq:4Dinnp} as $\R^4$, where we have for the orthogonal 2D planes the subspace bases\index{quaternions!subspace bases}:
\be 
  q_+\text{-basis: } \{f-g, 1+fg\},\qquad
  q_-\text{-basis: } \{f+g, 1-fg\}.
  \label{eq:q-q+basisfg}
\ee  
We can therefore use the following representation for every $q \in \H$ by means of four real coefficients\index{quaternions!orthogonal 2D planes split!real coefficients} $q_1, q_2, q_3, q_4 \in \R$
\begin{equation} 
  q = q_1 (1+fg) + q_2 (f-g) + q_3 (1-fg) + q_4 (f+g),
\end{equation}
where
\begin{align}
  q_1 &= \Scalar{q\inverse{(1+fg)}}, \qquad 
  q_2 = \Scalar{q\inverse{(f-g)}},\nonumber\\
  q_3 &= \Scalar{q\inverse{(1-fg)}}, \qquad 
  q_4 = \Scalar{q\inverse{(f+g)}}.
\end{align}
As an example we have for $f=\i, g=\j$ according to \eqref{eq:qpm}
\change{the coefficients for the decomposition with respect to the orthogonal basis \eqref{eq:fgsplitbasis}}
\begin{equation}
  q_1 = \frac{1}{2}(q_r+q_k), \quad
  q_2 = \frac{1}{2}(q_i-q_j), \quad
  q_3 = \frac{1}{2}(q_r-q_k), \quad
  q_4 = \frac{1}{2}(q_i+q_j).
\end{equation} 
Moreover, using 
\begin{equation}
  f-g = f (1+fg) = (1+fg)(-g),
  \quad 
  f+g = f(1-fg) = (1-fg)g,
  \label{eq:fgbasisrel}
\end{equation} 
we have the following left and right factoring properties
\begin{align} 
  q_+ 
  = q_1 (1+fg) + q_2 (f-g)
  &= (q_1+q_2f)(1+fg)
  \nonumber \\
  &= (1+fg)(q_1-q_2 g),
  \label{eq:q+} \\
  q_- 
  = q_3 (1-fg) + q_4 (f+g)
  &= (q_3+q_4f)(1-fg) 
  \nonumber \\
  &= (1-fg) (q_3+q_4 g).
  \label{eq:q-}
\end{align}

Equations \eqref{eq:fgplus} and \eqref{eq:fgminus} further show that the map $f\sandwich g$ rotates\index{quaternions!half-turn rotation} the $q_-$ plane by $180^{\circ}$ around the $q_+$ axis plane. We found that our interpretation of the map $f\sandwich g$ is in perfect agreement with Coxeter's notion of \emph{half-turn}\index{half-turn!Coxeter}\index{quaternions!half-turn!Coxeter} in \cite{Coxeter1946}. This opens the way for new types of QFTs, where the pair of square roots of $-1$ involved does not necessarily need to be orthogonal. 

Before suggesting a generalization of the QFT, we will establish a new set of very useful algebraic identities\index{quaternions!identites!algebraic}. Based on \eqref{eq:q+} and \eqref{eq:q-} we get for $\alpha,\beta \in \R$
\begin{gather} 
  e^{\alpha f} q e^{\beta g} 
  = e^{\alpha f} q_+ e^{\beta g} + e^{\alpha f} q_- e^{\beta g},
  \nonumber \\
  e^{\alpha f} q_+ e^{\beta g}
  = (q_1+q_2 f) e^{\alpha f} (1+fg) e^{\beta g}
  = e^{\alpha f} (1+fg) e^{\beta g} (q_1-q_2 g), 
  \nonumber \\
  e^{\alpha f} q_- e^{\beta g}
  = (q_3+q_4f) e^{\alpha f}(1-fg) e^{\beta g}
  = e^{\alpha f}(1-fg) e^{\beta g} (q_3+q_4g).
  \label{eq:eqpme}
\end{gather}
Using \eqref{eq:q+} again we obtain
\begin{align}
  e^{\alpha f} (1+fg) 
  &\stackrel{\phantom{\eqref{eq:q+}}}{=} (\cos \alpha + f \sin \alpha ) (1+fg) 
  \nonumber \\ 
  &\stackrel{\eqref{eq:q+}}{=} (1+fg) (\cos \alpha - g \sin \alpha )
  = (1+fg) e^{-\alpha g},
  \label{eq:eafleft}
\end{align}
where we set $q_1 = \cos \alpha$, $q_2 = \sin \alpha$ for applying \eqref{eq:q+}. Replacing in \eqref{eq:eafleft} $-\alpha \rightarrow \beta$ we get
\begin{equation}
  e^{-\beta f} (1+fg) 
  = (1+fg) e^{\beta g},
  \label{eq:ebfebg}
\end{equation} 
Furthermore, replacing in \eqref{eq:eafleft} $g \rightarrow -g$ and subsequently
$\alpha \rightarrow \beta$ we get
\begin{equation}
\change{
\begin{aligned}
  e^{\alpha f} (1-fg) &= (1-fg) e^{\alpha g},\\
  e^{\beta  f} (1-fg) &= (1-fg) e^{\beta  g},
\end{aligned}
}
  \label{eq:eafebfleft}
\end{equation} 
Applying \eqref{eq:q+}, \eqref{eq:eqpme}, \eqref{eq:eafleft} and  \eqref{eq:ebfebg} we can rewrite
\begin{align}
  e^{\alpha f} q_+ e^{\beta g}
  &\stackrel{\eqref{eq:eqpme}}{=} (q_1+q_2 f) e^{\alpha f} (1+fg) e^{\beta g}
  \stackrel{\eqref{eq:eafleft}}{=} (q_1+q_2 f) (1+fg) e^{(\beta-\alpha) g}\nonumber\\
  &\stackrel{\eqref{eq:q+}}{=} q_+ e^{(\beta-\alpha) g},
  \label{eq:eq+etoright}
\end{align} 
or equivalently as
\begin{align}
  e^{\alpha f} q_+ e^{\beta g}
  &\stackrel{\eqref{eq:eqpme}}{=} e^{\alpha f} (1+fg) e^{\beta g} (q_1-q_2 g)
  \stackrel{\eqref{eq:ebfebg}}{=} e^{(\alpha-\beta) f} (1+fg) (q_1-q_2 g)\nonumber\\
  &\stackrel{\eqref{eq:q+}}{=} e^{(\alpha-\beta) f} q_+ .
  \label{eq:eq+etoleft}
\end{align} 
In the same way by changing $g\rightarrow -g, \beta \rightarrow -\beta$ in \eqref{eq:eq+etoright} and \eqref{eq:eq+etoleft} we can rewrite
\begin{equation}
  e^{\alpha f} q_- e^{\beta g} 
  = e^{(\alpha+\beta) f} q_- 
  = q_- e^{(\alpha+\beta) g}.
\end{equation} 
The result is therefore\index{quaternions!OPS!exponential factors!identities}
\begin{equation}
  e^{\alpha f} q_{\pm} e^{\beta g} 
  = q_{\pm} e^{(\beta\mp\alpha) g}
  = e^{(\alpha\mp\beta) f}q_{\pm}.
  \label{eq:eqpmeres}
\end{equation}

\subsection{Orthogonal 2D Planes Split using One Pure Unit Quaternion\index{quaternions!OPS!single pure unit quaternion}}

We now treat the case for $g=f, f^2=-1$. We then have the map $f\sandwich f$, and the OPS split with respect to $f\in \H, f^2=-1$,
\be 
  q_{\pm} = \frac{1}{2} (q \pm f q f) .
  \label{eq:opsfqf}
\ee 
The pure quaternion $\i$ can be rotated by $R=\i(\i+f)$, see \eqref{eq:rotif}, into the quaternion unit $f$ and back. Therefore studying the map $\i\sandwich\i$ is\change{,}
up to the constant rotation between $\i$ and $f$\change{,} 
the same as studying $f\sandwich f$. This gives
\begin{equation}
  \i q \i 
  = \i (q_r + q_i \i + q_j \j + q_k \k)\i
  = -q_r - q_i \i + q_j \j + q_k \k .
  \label{eq:iqi}
\end{equation} 
The OPS with respect to $f=g=\i$ gives
\begin{equation}
  q_{\pm} = \frac{1}{2}(q\pm \i q \i), \quad 
  q_+ = q_j \j + q_k \k = (q_j + q_k \i)\j, \quad 
  q_- = q_r + q_i \i ,
  \label{eq:opsiqi}
\end{equation}  
where the $q_+$ plane is two-dimensional and manifestly orthogonal to the 2D $q_-$ plane. 
This form \eqref{eq:opsiqi} of the OPS is therefore identical to the quaternionic simplex/perplex split\index{quaternions!simplex and perplex split} of \cite{10.1109/TIP.2006.884955}.

For $g=f$ the $q_-$ plane is always spanned by $\{1,f\}$. The rotation operator $R=\i(\i+f)$,
with \change{squared norm} $|R|^2 = |\i(\i+f)|^2 = |(\i+f)|^2 = -(\i+f)^2$, rotates $\i$ into $f$ according to
\begin{gather} 
  R^{-1} \i R
  = \frac{\qconjugate{R}}{|R|^2}\i R
  = \frac{(\i+f)\i \i \i(\i+f)}{-(\i+f)^2}
  = \frac{(\i+f)\i (\i(-f)+1)f}{(\i+f)^2}
  \nonumber \\
  = \frac{(\i+f) (f+\i)f}{(\i+f)^2}
  = f . 
  \label{eq:rotif}
\end{gather} 
The rotation $R$ leaves $1$ invariant and thus rotates the whole $\{1,\i\}$ plane into the $q_-$ plane spanned by $\{1,f\}$. Consequently $R$ also rotates the $\{\j,\k\}$ plane into the $q_+$ plane spanned by $\{\j'= R^{-1} \j R, \,\k'=R^{-1} \k R\}$. We thus constructively obtain the fully \emph{orthonormal} 4D basis\index{quaternions!orthonormal basis} of $\H$ as
\be 
  \{1, f, \j', \k'\}
  = R^{-1}\{1,\i, \j, \k\}R, 
  \qquad R=\i(\i+f) ,
  \label{eq:onbasisff}
\ee 
for any chosen pure unit quaternion $f$. We further have\change{,} 
for the orthogonal 2D planes created in \eqref{eq:opsfqf} the subspace bases\index{quaternions!OPS!subspace bases}:
\be 
  q_+\text{-basis: } \{\j', \k'\},\qquad
  q_-\text{-basis: } \{1, f\}.
  \label{eq:q-q+basisff}
\ee  

The rotation $R$ (an orthogonal transformation!) of \eqref{eq:rotif} preserves the fundamental quaternionic orthonormality and the anticommutation relations
\be 
  f \j' = \k' = -\j' f ,
  \qquad
  \k' f = \j' = -f \k' 
  \qquad
  \j'\k' = f = -\k' \j' . 
  \label{eq:fq+acomm}
\ee 
Hence  
\begin{equation}
  f q f
  = f (q_+ + q_-) f = q_+ - q_- , 
  \,\,\,\mbox{ \ie }\,\,\,
  f q_{\pm} f = \pm q_{\pm} ,
\end{equation} 
which represents again a half-turn\index{quaternions!half-turn} by $180^{\circ}$ in the 2D $q_-$ plane around the 2D $q_+$ plane (as axis).  

\begin{figure}[p!]
\subfigure[][$q_+$ component]{\includegraphics[width=0.9\textwidth]{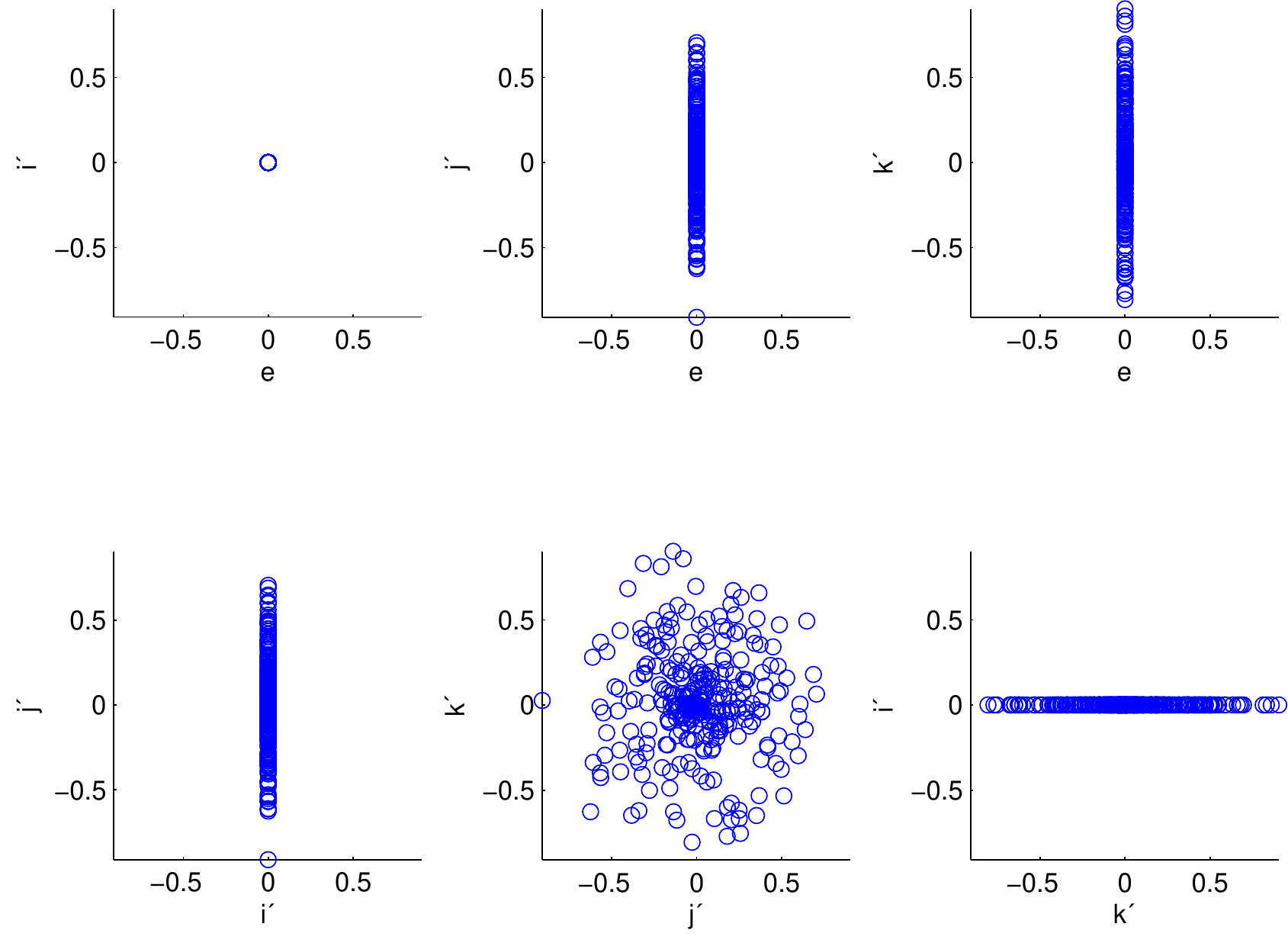}}
\subfigure[][$q_-$ component]{\includegraphics[width=0.9\textwidth]{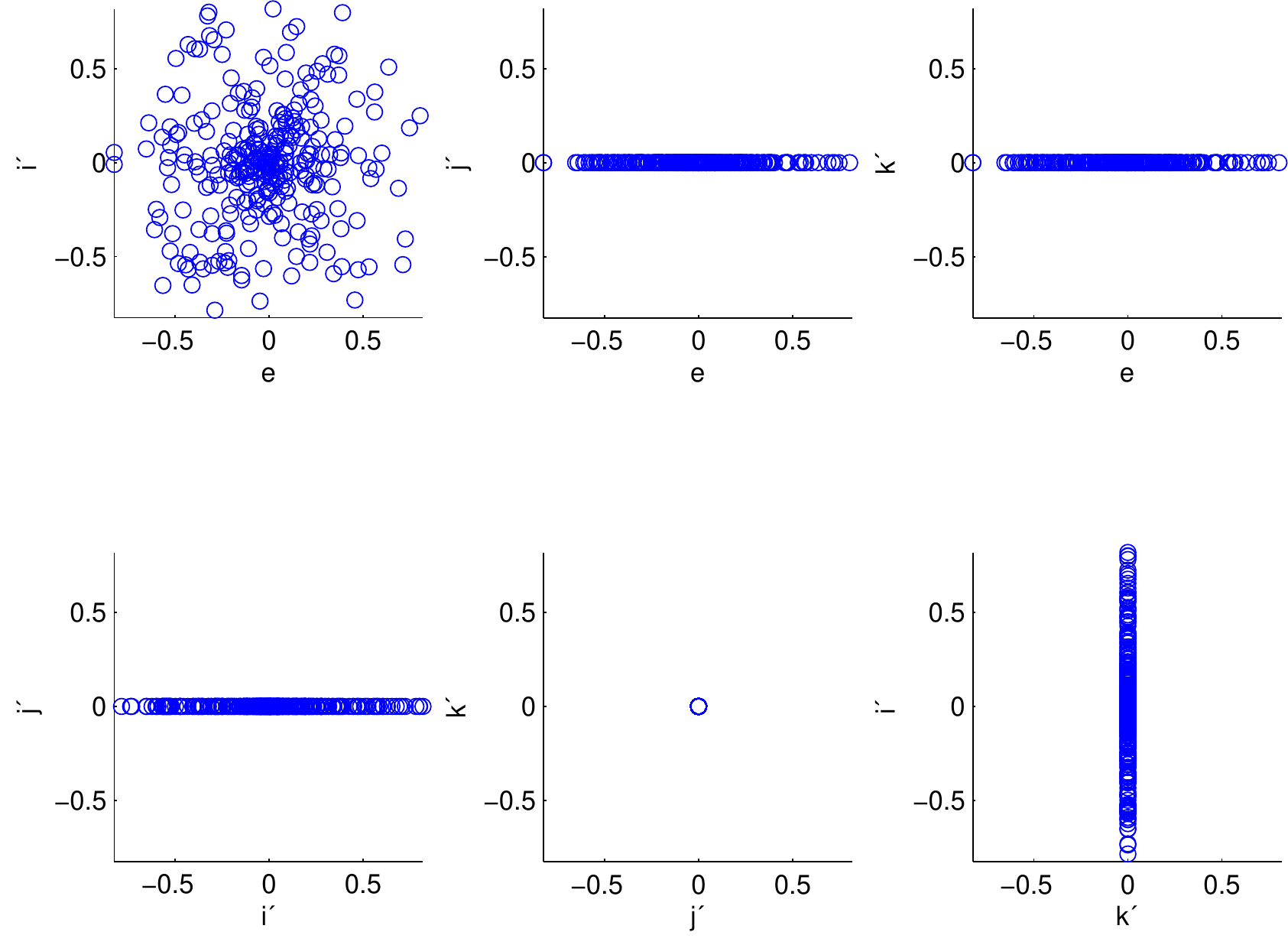}}
\caption{\label{fig:scatter4ff}4D scatter plot of quaternions decomposed using the orthogonal planes split
of \eqref{eq:opsfqf} {with one unit pure quaternion $f=\i'=\frac{1}{\sqrt{3}}( \i+\j+\k ) = g$}.}
\end{figure}

\begin{figure}
\includegraphics[width=0.6\textwidth]{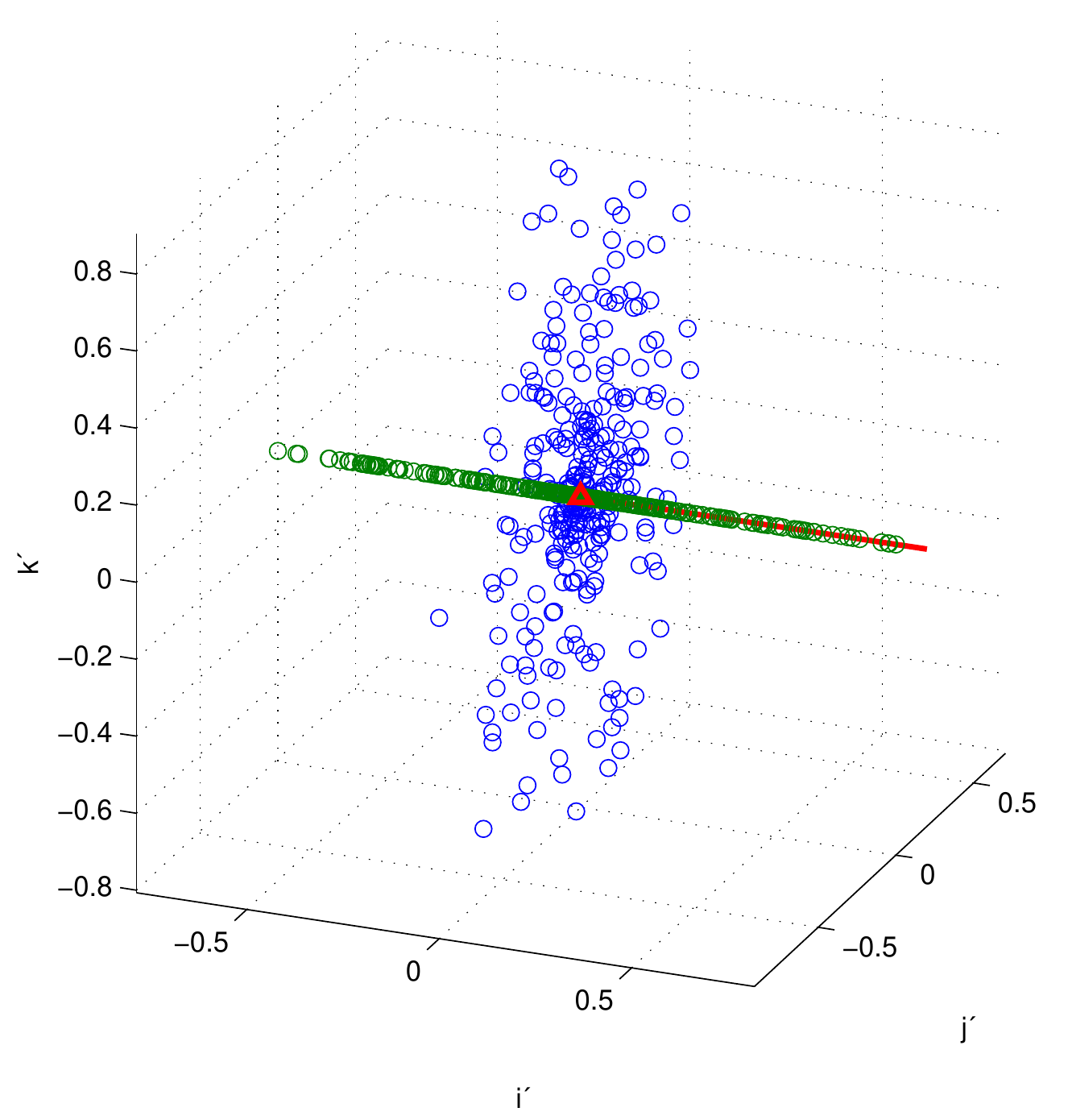}
\caption{\label{fig:scatter3ff}Scatter plot of vector parts of quaternions decomposed using the orthogonal planes split of \eqref{eq:opsfqf} {with one pure unit quaternion
$f=\i'=\frac{1}{\sqrt{3}}(\i+\j+\k)=g$}. The red line corresponds to the direction of $f$.}
\end{figure}

Figures \ref{fig:scatter4ff} and \ref{fig:scatter3ff} illustrate this decomposition for the case
where $f$ is the unit pure quaternion $\frac{1}{\sqrt{3}}(\i+\j+\k)$.
This decomposition corresponds (for pure quaternions) to the classical \emph{luminance-chrominance}\index{quaternions!luminance}\index{quaternions!chrominance}\index{luminance}\index{chrominance}\index{quaternions!!decomposition!luminance and chrominance}
decomposition used in colour image processing\index{quaternions!colour image!processing}\index{colour image!processing}, as illustrated, for example, in \cite[Figure 2]{10.1109/TIP.2006.884955}.
Three hundred unit quaternions randomly oriented in 4-space were decomposed.
Figure \ref{fig:scatter4ff} shows the {three hundred points} in 4-space,
{projected onto} the six orthogonal {planes}
{
$\{e,\i'\}, \{e,\j'\}, \{e,\k'\}, \{\i',\j'\}, \{\j',\k'\}, \{\k',\i'\}$
where $e=1$ and $\i'=f$, as given in \eqref{eq:onbasisff}.}
The six views on the left show the $q_+$ plane,
and the six on the right show the $q_-$ plane.

Figure \ref{fig:scatter3ff} shows the vector parts of the decomposed quaternions.
{The basis for the plot is $\{\i', \j', \k'\}$, where $\i'=f$ as given in \eqref{eq:onbasisff}.}
The green circles show
the components in the $\{1,f\}$ plane, which intersects the 3-space of the vector part only along the line
$f$ (which is the \emph{luminance} or \emph{grey line}\index{quaternions!grey line}\index{image processing!grey line}\index{colour image!grey line} of colour image pixels). The red line on the figure corresponds to
$f$. The blue circles show the components in the $\{\j^\prime,\k^\prime\}$ plane, which is entirely within
the 3-space. It is orthogonal to $f$ and corresponds to the \emph{chrominance plane}\index{quaternions!chrominance plane}\index{chrominance plane} of colour image processing.

The next question is the influence the current OPS \eqref{eq:opsfqf} has for left and right exponential factors of the form
\begin{equation}
  e^{\alpha f} q_{\pm}\,e^{\beta f} .
\end{equation} 
We learn from \eqref{eq:fq+acomm} that\index{quaternions!OPS!exponential factor!identity}
\begin{equation}
  e^{\alpha f} q_{\pm} e^{\beta f} 
  = e^{(\alpha\mp\beta)  f} q_{\pm}
  = q_{\pm} e^{(\beta \mp \alpha)  f},
\end{equation} 
which is identical to \eqref{eq:eqpmeres}, if we insert $g=f$ in \eqref{eq:eqpmeres}.  


Next, we consider $g=-f, f^2=-1$. We then have the map $f\sandwich (-f)$, and the OPS split with respect to $f,-f\in \H, f^2=-1$,
\be 
  q_{\pm} 
  = \frac{1}{2} (q \pm f q (-f)) 
  =\frac{1}{2} (q \mp f q f).
  \label{eq:opsfq-f}
\ee 
Again we can study $f=\i$ first, because for general pure unit quaternions $f$ the unit quaternion $\i$ can be rotated by \eqref{eq:rotif} into the quaternion unit $f$ and back. Therefore studying the map $\i\sandwich(-\i)$ is up to the constant rotation $R$ of \eqref{eq:rotif} the same as
studying $f\sandwich (-f)$. This gives the map
\begin{equation}
  \i q (-\i) 
  = \i (q_r + q_i \i + q_j \j + q_k \k)(-\i)
  = q_r + q_i \i - q_j \j - q_k \k .
  \label{eq:iq-i}
\end{equation} 
The OPS with respect to $f=\i,g=-\i$ gives
\begin{equation}
  q_{\pm} 
  = \frac{1}{2} (q \pm \i q (-\i)), \quad
  q_- = q_j \j + q_k \k = (q_j + q_k \i)\j, \quad 
  q_+ = q_r + q_i \i ,
  \label{eq:opsiq-i}
\end{equation}  
where, compared to $f=g=\i$, the 2D $q_+$ plane and the 2D $q_-$ plane\change{s} appear interchanged. 
The form \eqref{eq:opsiq-i} of the OPS is again identical to the quaternionic simplex/perplex\index{quaternions!simplex and perplex parts} split of \cite{10.1109/TIP.2006.884955}, but the simplex and perplex parts appear interchanged. 

For $g=-f$ the $q_+$ plane is always spanned by $\{1,f\}$. The rotation $R$ of \eqref{eq:rotif} rotates $\i$ into $f$ and leaves  $1$ invariant and thus rotates the whole $\{1,\i\}$ plane into the $q_+$ plane spanned by $\{1,f\}$. Consequently, $R$ of \eqref{eq:rotif} also rotates the $\{\j,\k\}$ plane into the $q_-$ plane spanned by $\{\j'= R^{-1} \j R, \,\k'=R^{-1} \k R\}$. 

We therefore have for the orthogonal 2D planes created in \eqref{eq:opsfq-f} the subspace bases\index{quaternions!OPS!subspace bases}:
\be 
  q_+\text{-basis: } \{1, f\},\qquad
  q_-\text{-basis: } \{\j', \k'\}.
  \label{eq:q-q+basisf-f}
\ee  

We again obtain the fully \emph{orthonormal} 4D basis\index{quaternions!orthonormal basis} \eqref{eq:onbasisff} of $\H$, preserving the fundamental quaternionic orthonormality and the anticommutation relations \eqref{eq:fq+acomm}.

Hence for \eqref{eq:opsfq-f}
\begin{equation}
  f q (-f)
  = f (q_+ + q_-) (-f) = q_+ - q_- , 
  \,\,\,\mbox{ \ie }\,\,\,
  f q_{\pm} (-f) = \pm q_{\pm} ,
\end{equation} 
which represents again a half-turn\index{quaternions!half-turn} by $180^{\circ}$ in the 2D $q_-$ plane around the 2D $q_+$ plane (as axis).  

The remaining question is the influence the current OPS \eqref{eq:opsfq-f} has for left and right exponential factors of the form
\begin{equation}
  e^{\alpha f} q_{\pm} e^{-\beta f} .
\end{equation} 
We learn from \eqref{eq:fq+acomm} that
\begin{equation}
  e^{\alpha f} q_{\pm} e^{-\beta f} 
  = e^{(\alpha\mp\beta)  f} q_{\pm}
  = q_{\pm} e^{-(\beta \mp \alpha) f},
\end{equation} 
which is identical to \eqref{eq:eqpmeres}, if we insert $g=-f$ in \eqref{eq:eqpmeres}. 

For \eqref{eq:eqpmeres} therefore, we do not any longer need to distinguish the cases $f\neq \pm g$ and $f=\pm g$. 
This motivates us to a general OPS definition for any pair of pure quaternions $f,g$, and we get a general lemma. 

\begin{defn}[General orthogonal 2D planes split\index{quaternions!orthogonal 2D planes split!general}]
\label{df:genOPS}
Let $f,g \in \H$ be an arbitrary pair of pure quaternions $f,g$, $f^2=g^2=-1$, including the cases $f = \pm g$. The general OPS is then defined with respect to the two pure unit quaternions $f, g$ as 
\begin{equation}
  \label{eq:genOPS}
  q_{\pm} = \frac{1}{2}(q \pm f q g).
\end{equation} 
\end{defn}

\begin{rem}
The three generalized OPS \eqref{eq:opsfqg}, \eqref{eq:opsfqf}, and \eqref{eq:opsfq-f} are formally identical and are now subsumed in  \eqref{eq:genOPS} of Definition \ref{df:genOPS}, where the values $g=\pm f$ are explicitly included, \ie any pair of pure unit quaternions $f,g\in \H$, $f^2=g^2=-1$, is admissible.
\end{rem}

\begin{lem}\label{lm:expqexp}
With respect to the general OPS of Definition \ref{df:genOPS} we have for left and right exponential factors the identity\index{quaternions!exponential factors!identity}
\begin{equation}
  e^{\alpha f} q_{\pm} e^{\beta g} 
  = q_{\pm} e^{(\beta\mp\alpha) g}
  = e^{(\alpha\mp\beta) f}q_{\pm}.
  \label{eq:lemexpqexp}
\end{equation} 
\end{lem}

\subsection{Geometric Interpretation of Left and Right Exponential Factors in $f$, $g$ \label{ssc:geointerp}}

We obtain the following general \emph{geometric interpretation}\index{quaternions!OPS!geometric interpretation}. The map $f\sandwich g$ always represents a rotation by angle $\pi$ in the $q_-$ plane (around the $q_+$ plane), the map $f^t\sandwich g^t$, $t\in \R$, similarly represents a rotation by angle $t\pi$ in the $q_-$ plane (around the $q_+$ plane as axis).
Replacing\footnote{\change{Alternatively and equivalently we could replace
$f\rightarrow -f$ instead of $g \rightarrow -g$.}}
$g\rightarrow -g$ in the map $f\sandwich g$ we further find that 
\begin{equation}
  f q_{\pm} (-g) = \mp q_{\pm}.
\end{equation} 
Therefore the map $f\sandwich(-g) = f\sandwich\inverse{g}$, because $\inverse{g}=-g$, represents a rotation by angle $\pi$ in the $q_+$ plane (around the $q_-$ plane), exchanging the roles of 2D rotation plane and 2D rotation axis. Similarly, the map $f^s\sandwich g^{-s}$, $s\in \R$, represents a rotation by angle $s\pi$ in the $q_+$ plane (around the $q_-$ plane as axis). 

The product of these two rotations gives
\begin{gather}
  f^{t+s} q g^{t-s} 
  = e^{(t+s)\frac{\pi}{2}f}qe^{(t-s)\frac{\pi}{2}g} 
  = e^{{\alpha}f}qe^{{\beta}g} , 
  \nonumber \\
  \alpha = (t+s)\frac{\pi}{2}, \qquad 
  \beta = (t-s)\frac{\pi}{2} ,
  \label{eq:rotab}
\end{gather} 
where based on \change{\eqref{eq:unit}}
we used the identities $f=e^{\frac{\pi}{2}f}$ and $g=e^{\frac{\pi}{2}g}$.

The \emph{geometric interpretation}\index{quaternions!OPS!geometric interpretation} of \eqref{eq:rotab} is a rotation by angle $\alpha+\beta$ in the $q_-$ plane (around the $q_+$ plane), and a second rotation\index{quaternions!OPS!geometric interpretation!rotation} by angle $\alpha-\beta$ in the $q_+$ plane (around the $q_-$ plane). For $\alpha=\beta=\pi/2$ we recover the map $f\sandwich g$, and for $\alpha=-\beta=\pi/2$ we recover the map $f\sandwich\inverse{g}$.

\subsection{Determination of $f,g$ for Given Steerable Pair of Orthogonal 2D planes\index{quaternions!OPS!determination from given planes}}
\label{sec:determinesteer}

Equations \eqref{eq:q-q+basisfg}, \eqref{eq:q-q+basisff}, and \eqref{eq:q-q+basisf-f} tell us how the pair of pure unit quaternions $f,g \in \H$ used in the general OPS of Definition \ref{df:genOPS}, leads to an explicit basis for the resulting two orthogonal 2D planes, the $q_+$ plane and the $q_-$ plane. We now ask the \emph{opposite} question: how can we determine from a given steerable pair of orthogonal 2D planes in $\H$ the pair of pure unit quaternions $f,g \in \H$, which splits $\H$ exactly into this given pair of orthogonal 2D planes?

To answer this question, we first observe that in a 4D space it is sufficient to know only one 2D plane explicitly, specified \eg by a pair of orthogonal unit quaternions $a,b \in \H$, $|a|=|b|=1$, and without restriction of generality $b^2=-1$, \ie $b$ can be a pure unit quaternion {$b=\bmu(b)$}. But for $a=\cos \alpha + {\bmu(a)}\sin \alpha$, compare \eqref{eq:unitq}, we must distinguish $\Scalar{a} = \cos \alpha \neq 0$ and $\Scalar{a} = \cos \alpha = 0$, \ie of $a$ also being a pure quaternion with $a^2=-1$. The second orthogonal 2D plane is then simply the \emph{orthogonal complement} in $\H$ to the $a,b$ plane.  

Let us first treat the case $\Scalar{a} = \cos \alpha \neq 0$. We set 
\be 
  f := ab, \qquad
  g := \qconjugate{a}b.
  \label{eq:detfgSbnot0}
\ee 
With this setting we get for the basis of the $q_-$ plane
\begin{align} 
  f+g &= ab+\qconjugate{a}b = 2 \Scalar{a} b, 
  \nonumber \\
  1-fg 
  &= 1-ab\qconjugate{a}b 
  = 1-a^2b^2 = 1+a^2 
  \nonumber \\
  &= 1+\cos^2 \alpha-\sin^2 \alpha + 2 {\bmu(a)} \cos \alpha \sin \alpha 
  \nonumber \\
  &= 2 \cos \alpha (\cos \alpha + {\bmu(a)} \sin \alpha)
  = 2 \Scalar{a} a.
  \label{eq:q-detfgSbnot0}
\end{align} 
For the equality 
  $ab\qconjugate{a}b = a^2b^2$
we used the orthogonality of $a,b$, which means that the vector part of $a$ must be orthogonal to the pure unit quaternion $b$, \ie it must anticommute with $b$
\be 
  ab = b\qconjugate{a}, \qquad 
  ba = \qconjugate{a}b .
  \label{eq:abcomm}
\ee 
Comparing \eqref{eq:q-q+basisfg} and \eqref{eq:q-detfgSbnot0}, the plane spanned by the two orthogonal unit quaternions $a,b \in \H$ is indeed the $q_-$ plane for $\Scalar{a} = \cos \alpha \neq 0$. The orthogonal $q_+$ plane is simply given by its basis vectors \eqref{eq:q-q+basisfg}, inserting \eqref{eq:detfgSbnot0}. This leads to the pair of orthogonal unit quaternions $c,d$ for the $q_+$ plane as
\be
  c = \frac{f-g}{|f-g|}
  = \frac{ab - \qconjugate{a}b}{|(a - \qconjugate{a})b|}
  = \frac{a - \qconjugate{a}}{|a - \qconjugate{a}|}b
  = {\bmu(a)} b,
\ee 
\be
  d = \frac{1+fg}{|1+fg|}
  = \frac{f-g}{|f-g|}g
  = cg 
  = {\bmu(a)} bg
  = {\bmu(a)} b\qconjugate{a}b
  = {\bmu(a)} ab^2
  = -{\bmu(a)} a ,
  \label{eq:dcalc}
\ee 
where we have used \eqref{eq:fgbasisrel} for the second, and \eqref{eq:abcomm} for the sixth equality in \eqref{eq:dcalc}.

Let us also verify that $f,g$ of \eqref{eq:detfgSbnot0} are both pure unit quaternions using \eqref{eq:abcomm}
\be 
  f^2 = abab = (a\qconjugate{a})bb = -1,
  \qquad 
  g^2 = \qconjugate{a}b\qconjugate{a}b = (\qconjugate{a}a)bb = -1. 
\ee

Note, that if we would set in \eqref{eq:detfgSbnot0} for $g := -\qconjugate{a}b$, then the $a,b$ plane would have become the $q_+$ plane instead of the $q_-$ plane. We can therefore determine by the sign in the definition of $g$, which of the two OPS planes the $a,b$ plane is to represent.

For both $a$ and $b$ being two pure orthogonal quaternions, we can again set 
\be 
  f := ab \Rightarrow f^2 = abab = -a^2b^2 = -1, \qquad
  g := \qconjugate{a}b = -ab = -f,
\ee 
where due to the orthogonality of the pure unit quaternions $a,b$ we were able to use $ba = -ab$. In this case $f=ab$ is thus also shown to be a pure unit quaternion. Now the $q_-$ plane of the corresponding OPS \eqref{eq:opsfq-f} is spanned by $\{a,b\}$, whereas the $q_+$ plane is spanned by $\{1,f\}$. Setting instead $g:=-\qconjugate{a}b = ab = f$, the $q_-$ plane of the corresponding OPS \eqref{eq:opsfqf} is spanned by $\{1,f\}$, wheras the $q_+$ plane is spanned by $\{a,b\}$. 

We summarize our results in \change{the following Theorem.}\index{quaternions!OPS!determination from given plane}
\begin{thm}
\textbf{(Determination of \ensuremath{f,g} from given steerable  2D plane)}\label{th:abdetfg}
Given any 2D plane in $\H$ in terms of two unit quaternions $a,b$, where $b$ is without restriction of generality pure, \ie $b^2=-1$, we can make the given plane the $q_-$ plane of the OPS $q_{\pm}=\frac{1}{2}(q\pm fqg)$, by setting
\be 
  f := ab, \qquad
  g := \qconjugate{a}b.
  \label{eq:theoremdetfg}
\ee 
For $\Scalar{a}\neq 0$
the orthogonal $q_+$ plane is fully determined by the orthogonal unit quaternions
\be 
  c =  {\bmu(a)}b,
  \qquad
  d = -{\bmu(a)}a .
  \label{eq:cdindetfg}
\ee 
{where $\bmu(a)$ is as defined in} \eqref{eq:unitq}.
For $\Scalar{a}= 0$ the orthogonal $q_+$ plane with basis $\{1,f \}$ is instead fully determined by $f=-g=ab$.

Setting alternatively
\be 
  f := ab, \qquad
  g := -\qconjugate{a}b.
  \label{eq:theoremdetfg-}
\ee 
makes the given $a,b$ plane the $q_+$ plane instead. For $\Scalar{a}\neq 0$
the orthogonal $q_-$ plane is then fully determined by \eqref{eq:theoremdetfg-} and \eqref{eq:q-q+basisfg}, with the same orthogonal unit quaternions $c = \bmu(a)b, d = -\bmu(a)a$ as in \eqref{eq:cdindetfg}.
For $\Scalar{a}= 0$ the orthogonal $q_-$ plane with basis $\{1,f \}$ is then instead fully determined by $f=g=ab$.
\end{thm}

\begin{figure}[p!]
\subfigure[][$q_+$ component]{\includegraphics[width=0.9\textwidth]{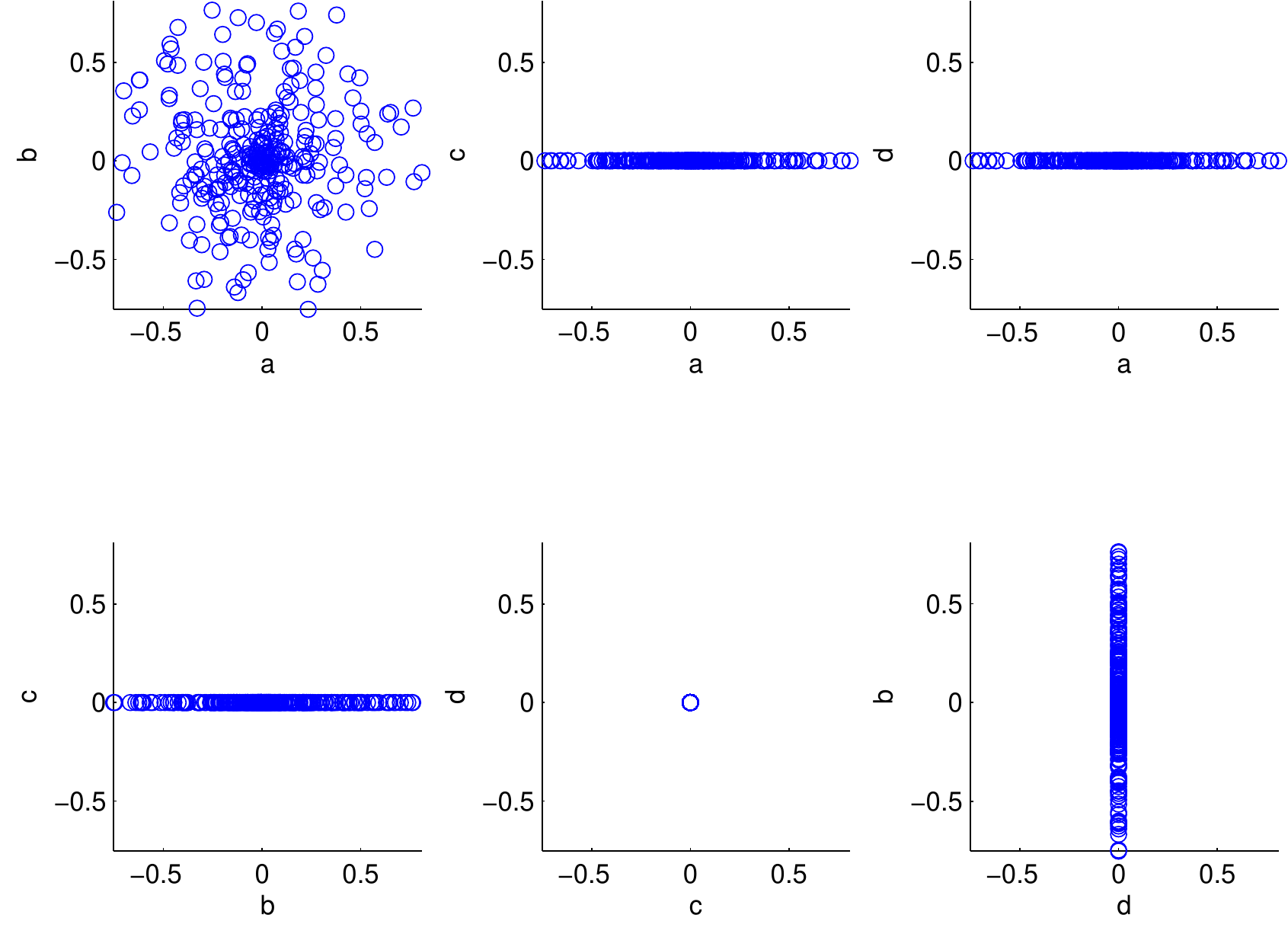}}
\subfigure[][$q_-$ component]{\includegraphics[width=0.9\textwidth]{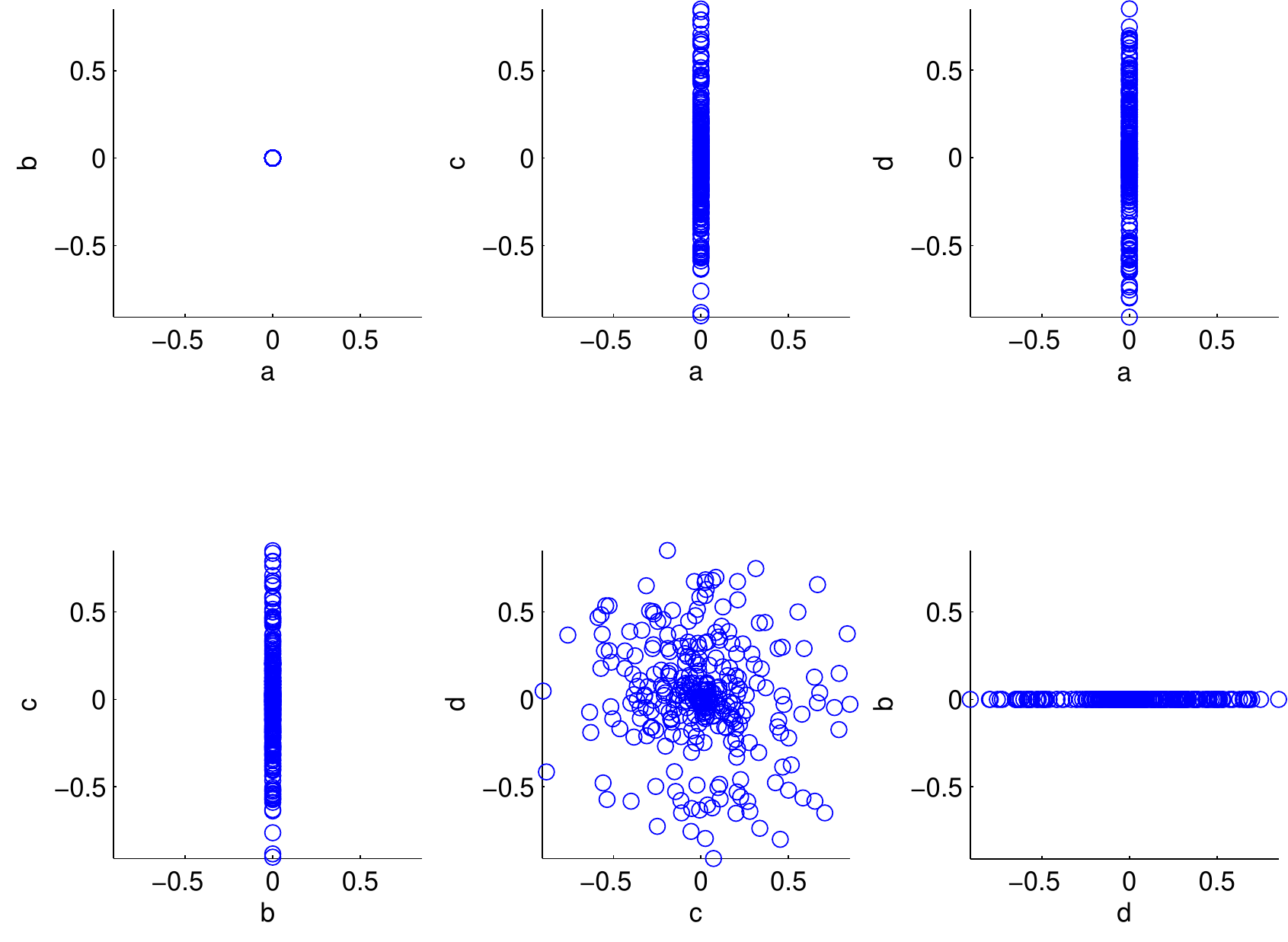}}
\caption{\label{fig:scatter4fg}4D scatter plot of quaternions decomposed using the orthogonal planes split
of Definition \ref{df:genOPS}.}
\end{figure}

\begin{figure}
\includegraphics[width=0.6\textwidth]{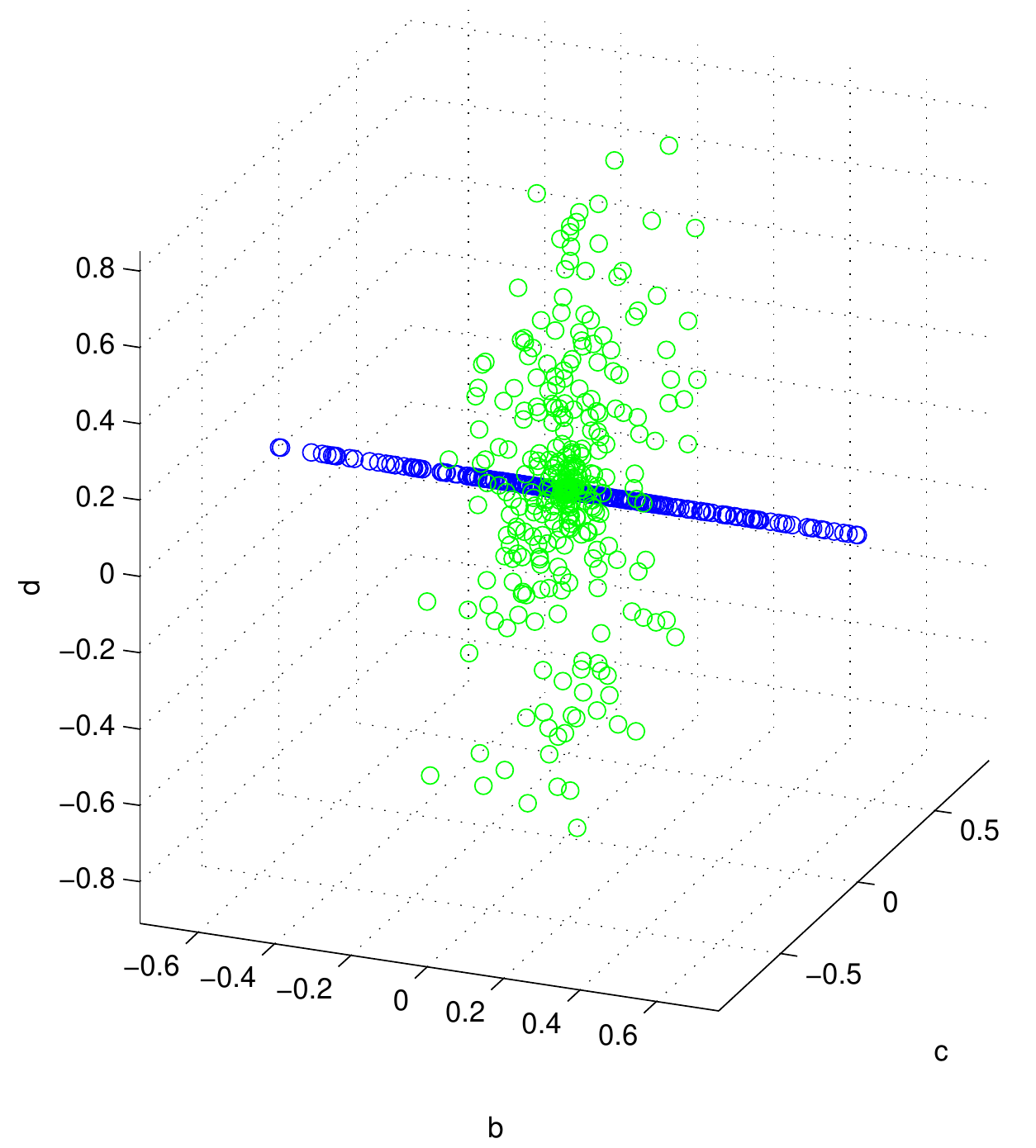}
\caption{\label{fig:scatter3fg}Scatter plot of vector parts of quaternions decomposed using the orthogonal
planes split of Definition \ref{df:genOPS}.}
\end{figure}

An illustration of the decomposition is given in \figurename{}s~ \ref{fig:scatter4fg} and \ref{fig:scatter3fg}.
Again, three hundred unit pure quaternions randomly oriented in 4-space have been
decomposed into two sets using the decomposition of Definition \ref{df:genOPS} and two unit pure quaternions
$f$ and $g$ computed as in Theorem \ref{th:abdetfg}.
{$b$ was the pure unit quaternion $\frac{1}{\sqrt{3}}(\i+\j+\k)$
and $a$ was the full unit quaternion $\frac{1}{\sqrt{2}} + \frac{1}{2}(\i-\j)$.}
%
{$c$ and $d$ were computed by \eqref{eq:cdindetfg} as $c=\bmu(a)b$ and $d=-\bmu(a)a$.}

Figure \ref{fig:scatter4fg} shows the three hundred points in 4-space,
projected onto the six orthogonal planes
$\{c,d\}, \{c,b\}, \{c,a\}, \{d,b\}, \{b,a\}, \{a,d\}$
where the orthonormal 4-space basis
$\{c, d, b, a\}
= \{(f-g)/|f-g|, (1+fg)/|1+fg|, (f+g)/|f+g|, (1-fg)/|1-fg|\}$.
The six views on the left show the $q_+$ plane,
and the six on the right show the $q_-$ plane.
Figure \ref{fig:scatter3fg} shows the vector parts of the decomposed quaternions.

\section{New QFT Forms: OPS-QFTs with Two Pure Unit Quaternions $f,g$\index{quaternions!QFT!new forms}}
\label{sc:NewQFT}

\subsection{Generalized OPS leads to new steerable type of QFT\label{ssc:steerQFT}}

We begin with a straightforward generalization\index{quaternions!QFT!generalization} of the (double sided form of the) QFT \cite{10.1007/s00006-007-0037-8,EH:DirUP_QFT} in $\H$ by replacing $\i$ \change{with} $f$ and $\j$ \change{with} $g$
defined as

\begin{defn}\textbf{(QFT with respect to two pure unit quaternions \ensuremath{f,g})}\label{df:QFTfg}\\
Let $f,g \in \H$, $f^2=g^2=-1$, be any two pure unit quaternions.\index{quaternions!QFT!two pure unit quaternions} The quaternion Fourier transform with respect to $f,g$ is 
\begin{equation}
  \label{eq:QFTfg}
  \mathcal{F}^{f,g}\{ h \}(\bomega)
  = \int_{\R^2} e^{-f x_1\omega_1} h(\bvect{x}) \,e^{-g x_2\omega_2} d^2\bvect{x},
\end{equation}
where $h\in L^1(\R^2,\H)$, $d^2\bvect{x} = dx_1dx_2$ and $\bvect{x}, \bomega \in \R^2 $.
\end{defn}
Note, that the pure unit quaternions $f,g$ in Definition \ref{df:QFTfg} do not need to be orthogonal, and that the cases $f=\pm g$ are fully included.

Linearity of the integral \eqref{eq:QFTfg} allows us to use the OPS split $h=h_- + h_+$
\begin{gather}
  \mathcal{F}^{f,g}\{ h \}(\bomega)
  = \mathcal{F}^{f,g}\{ h_- \}(\bomega) + \mathcal{F}^{f,g}\{ h_+ \}(\bomega)
  \nonumber \\
  = \mathcal{F}^{f,g}_-\{ h \}(\bomega) + \mathcal{F}_+^{f,g}\{ h \}(\bomega),
\end{gather} 
since by their construction the operators of the Fourier transformation $\mathcal{F}^{f,g}$, and of the OPS with respect to $f,g$ commute. 
From Lemma \ref{lm:expqexp} follows 

\begin{thm}[QFT of $h_{\pm}$]
\label{th:fpmtrafo}
The QFT of the $h_{\pm}$ OPS split parts\index{quaternions!QFT!split parts!transformation}, with respect to two unit quaternions $f,g$, of a quaternion module function 
$h \in L^1(\R^2,\H)$ have the quasi-complex forms\index{quaternions!QFT!quasi-complex}
\begin{gather}
  \mathcal{F}^{f,g}_{\pm}\{h\} 
  = \mathcal{F}^{f,g}\{h_{\pm}\} 
  \nonumber \\
  \stackrel{}{=} \int_{\R^2}
    h_{\pm}e^{-g (x_2\omega_2 \mp x_1\omega_1)}d^2x
  \stackrel{}{=} \int_{\R^2}
    e^{-f (x_1\omega_1 \mp x_2\omega_2)}h_{\pm}d^2x \,\, .
\end{gather}
\end{thm} 

\begin{rem}
The quasi-complex forms in Theorem \ref{th:fpmtrafo} allow \change{us} to establish \emph{discretized}\index{quaternions!QFT!discrete} and \emph{fast}\index{quaternions!QFT!fast} versions of the QFT of Definition \ref{df:QFTfg} as sums of two complex discretized and fast Fourier transformations (FFT), respectively. 
\end{rem}

\begin{figure}
\begin{center}
\includegraphics[width=0.5\textwidth]{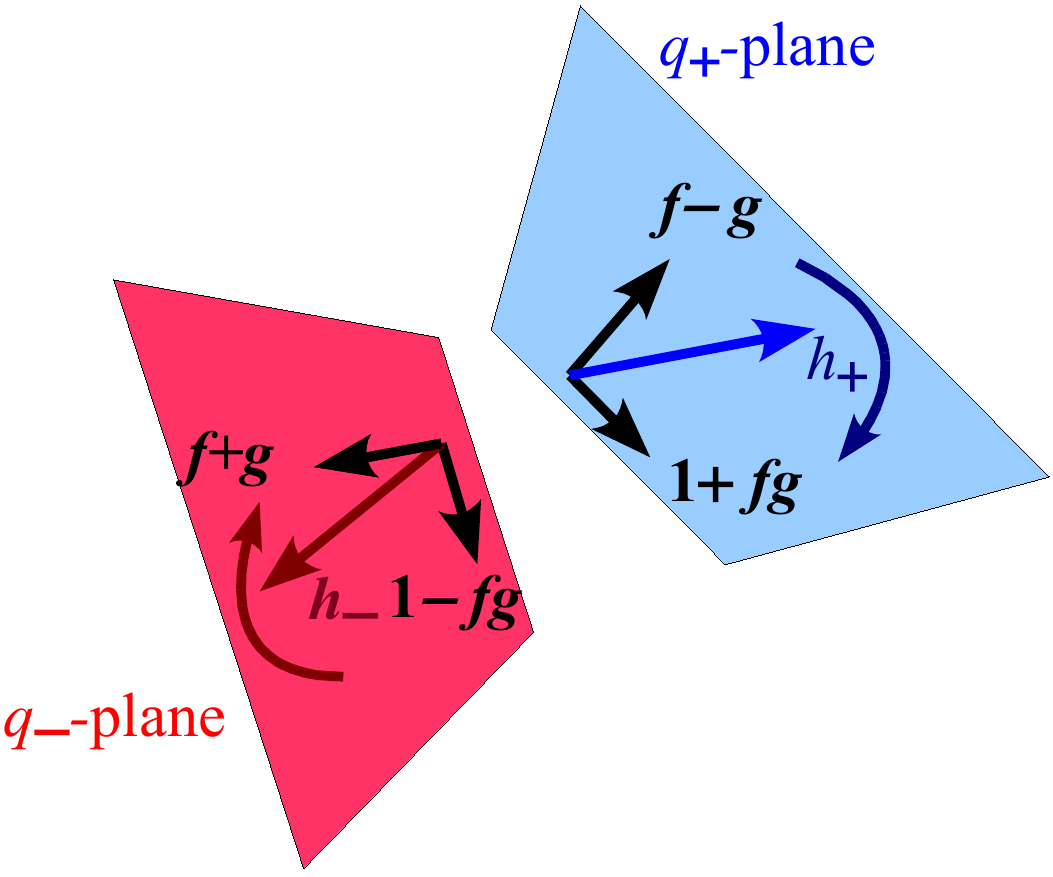}
\caption{\label{fig:orthoplanes} Geometric interpretation of integrand of QFT${}^{f,g}$ in Definition \ref{df:QFTfg} in terms of local phase rotations\index{quaternions!QFT!phase rotations!local} in $q_{\pm}$-planes.}
\end{center}
\end{figure}

We can now give a \emph{geometric interpretation}\index{quaternions!QFT!geometric interpretation} of the integrand \change{of} the QFT${}^{f,g}$ in Definition \ref{df:QFTfg} in terms of local phase rotations\index{quaternions!QFT!local phase rotations}, compare Section \ref{ssc:geointerp}. 
The integrand product 
\begin{equation}
  e^{-f x_1\omega_1} h(\bvect{x}) \,e^{-g x_2\omega_2}
\end{equation} 
\change{represents a \emph{local rotation} by the}
phase angle $-(x_1\omega_1+x_2\omega_2)$ in the $q_-$ plane, and by the phase angle\index{quaternions!QFT!phase angle}  $-(x_1\omega_1-x_2\omega_2) = x_2\omega_2-x_1\omega_1$ in the orthogonal $q_+$ plane, compare \figurename~\ref{fig:orthoplanes},
which depicts two completely orthogonal planes in four dimensions.


Based on Theorem \ref{th:abdetfg} the two phase rotation planes\index{quaternions!QFT!phase rotation planes} (analysis planes\index{quaternions!QFT!analysis planes}) can be freely \emph{steered} by defining the two pure unit quaternions $f,g$ used in Definition \ref{df:QFTfg} according to \eqref{eq:theoremdetfg} or \eqref{eq:theoremdetfg-}.

\subsection{Two Phase Angle Version of QFT\label{ssc:2pangleQFT}}

The above newly gained geometric understanding motivates us to propose a further new version of the QFT${}^{f,g}$, with a straightforward two \emph{phase angle} interpretation\change{.}

\begin{defn}\textbf{(Phase angle QFT\index{quaternions!QFT!phase angle transformation} with respect to  \ensuremath{f,g})}\label{df:QFT2angle}\\
Let $f,g \in \H$, $f^2=g^2=-1$, be any two pure unit quaternions. 
The phase angle quaternion Fourier transform with respect to $f,g$ is 
\begin{equation}
  \label{eq:QFT2angle}
  \mathcal{F}_D^{f,g}\{ h \}(\bomega)
  = \int_{\R^2} e^{-f \frac{1}{2}(x_1\omega_1+x_2\omega_2)} h(\bvect{x}) 
    \,e^{-g \frac{1}{2}(x_1\omega_1-x_2\omega_2)} d^2\bvect{x}.
\end{equation}
where again $h\in L^1(\R^2,\H)$, $d^2\bvect{x} = dx_1dx_2$ and $\bvect{x}, \bomega \in \R^2 $. 
\end{defn}

The \emph{geometric interpretation} of the integrand of \eqref{eq:QFT2angle} is a local phase rotation\index{quaternions!QFT!local phase rotation} by angle 
$-(x_1\omega_1+x_2\omega_2)/2- (x_1\omega_1-x_2\omega_2)/2 = - x_1\omega_1$ 
in the $q_-$ plane, and a second local phase rotation by angle 
$-(x_1\omega_1+x_2\omega_2)/2+ (x_1\omega_1-x_2\omega_2)/2 = - x_2\omega_2$
in the $q_+$ plane, compare Section \ref{ssc:geointerp}.

If we apply the OPS${}^{f,g}$ split to \eqref{eq:QFT2angle} we obtain the following theorem.

\begin{thm}[Phase angle QFT\index{quaternions!QFT!phase angle transform!split parts} of $h_{\pm}$]
\label{th:2anglefpmtrafo}
The phase angle QFT of Definition \ref{df:QFT2angle} applied to the $h_{\pm}$ OPS split parts, with respect to two pure unit quaternions $f,g$, of a quaternion module function $h \in L^1(\R^2,\H)$ leads to the quasi-complex expressions
\begin{equation}
  \mathcal{F}^{f,g}_{D+}\{h\} 
  = \mathcal{F}^{f,g}_{D}\{h_+\}
  \stackrel{}{=} \int_{\R^2}
    h_{+}e^{+g x_2\omega_2 }d^2x
  \stackrel{}{=} \int_{\R^2}
    e^{-f x_2\omega_2}h_{+}d^2x \,\, ,
\end{equation}
\begin{equation}
  \mathcal{F}^{f,g}_{D-}\{h\} 
  =\mathcal{F}^{f,g}_{D}\{h_-\} 
  \stackrel{}{=} \int_{\R^2}
    h_{-}e^{-g x_1\omega_1 }d^2x
  \stackrel{}{=} \int_{\R^2}
    e^{-f x_1\omega_1}h_{-}d^2x \,\, ,
\end{equation}
\end{thm} 

Note that based on Theorem \ref{th:abdetfg} the two phase rotation planes (analysis planes) are again freely steerable\index{quaternions!QFT!phase angle transform!steerable}. 

Theorem \ref{th:2anglefpmtrafo} allows \change{us} to establish \emph{discretized}\index{quaternions!QFT!phase angle transform!discrete} and \emph{fast}\index{quaternions!QFT!phase angle transform!fast} versions of the phase angle QFT of Definition \ref{df:QFT2angle} as sums of two complex discretized and fast Fourier transformations (FFT), respectively. 

The maps $f\sandwich g$ considered so far did not involve \emph{quaternion conjugation} $q\rightarrow \qconjugate{q}$. In the following we investigate maps which additionally 
conjugate the argument, 
\ie of type $f\qconjugate{\sandwich}g$, which are also \emph{involutions}.

\section{Involutions and QFTs Involving Quaternion Conjugation\label{sc:conjQFT}\index{quaternions!QFT!quaternion conjugation}}

\subsection{Involutions Involving Quaternion Conjugations\index{quaternions!involution!quaternion conjugation}}

The simplest case is quaternion conjugation itself
\begin{equation}
  q \rightarrow \qconjugate{q} = q_r - q_i \i -q_j \j - q_k \k ,
\end{equation} 
which can be interpreted as a \emph{reflection at the real line}\index{quaternions!line reflection!real line} $q_r$. The real line through the origin remains pointwise invariant, while every other point in the 3D hyperplane of pure quaternions is reflected to the opposite side of the real line. The related involution
\begin{equation}
  q \rightarrow -\qconjugate{q} = - q_r + q_i \i +q_j \j + q_k \k ,
  \label{eq:3Dhypref}
\end{equation} 
is the \emph{reflection at the 3D hyperplane}\index{quaternions!reflection!hyperplane} of pure quaternions (which stay invariant), \ie only the real line is changed into its negative $q_r \rightarrow -q_r$.

Similarly any pure unit quaternion factor like $\i$ in the map
\begin{equation}
  q \rightarrow \i\,\qconjugate{q}\i = -q_r + q_i \i -q_j \j - q_k \k ,
\end{equation} 
leads to a reflection at the (pointwise invariant) {line through the origin with direction $\i$,} while the map 
\begin{equation}
  q \rightarrow -\i\,\qconjugate{q}\i = q_r - q_i \i +q_j \j + q_k \k ,
\end{equation} 
leads to a reflection at the invariant 3D hyperplane orthogonal\index{quaternions!reflection!hyperplane!invariant} to the {line through the origin with direction $\i$.}
\change{The map}
\begin{equation}
  q \rightarrow f\,\qconjugate{q}f,
\end{equation} 
leads to a reflection at the (pointwise invariant) line with direction $f$ through the origin, while the map 
\begin{equation}
  q \rightarrow -f\,\qconjugate{q}f ,
\end{equation} 
leads to a reflection at the invariant 3D hyperplane orthogonal to the line with direction $f$ through the origin.

Next we turn to a map of the type 
\begin{equation}
  q \rightarrow - e^{\alpha f}\qconjugate{q}e^{\alpha f} .
\end{equation} 
Its set of pointwise invariants is given by 
\begin{align}
  q = - e^{\alpha f}\qconjugate{q}e^{\alpha f} 
  &\Leftrightarrow \,\,\,
  e^{-\alpha f}q = - \qconjugate{q}e^{\alpha f}\,\,\,
  \Leftrightarrow \,\,\,
  e^{-\alpha f}q + \qconjugate{q}e^{\alpha f} = 0
  \nonumber\\
  &\Leftrightarrow \,\,\,
  \Scalar{\qconjugate{q}e^{\alpha f}} = 0\,\,\,
  \Leftrightarrow \,\,\,
  q \perp e^{\alpha f}.
\end{align} 
We further observe that
\begin{equation}
  e^{\alpha f} \rightarrow - e^{\alpha f}e^{-\alpha f}e^{\alpha f} 
  = - e^{\alpha f}.
\end{equation} 
The map $-a\qconjugate{\sandwich}a$, with unit quaternion $a=e^{\alpha f}$, therefore represents a reflection at the invariant 3D hyperplane orthogonal to the line through the origin with direction $a$. 

Similarly, the map $a\qconjugate{\sandwich}a$, with unit quaternion $a=e^{\alpha f}$, then represents a reflection at the (pointwise invariant) line\index{quaternions!line reflection!invariant line!pointwise} with direction $a$ through the origin. 

The combination of two such reflections (both at 3D hyperplanes, or both at lines), given by unit quaternions $a,b$, leads  to a rotation\index{quaternions!rotation}
\begin{gather}
  -b\qconjugate{-a \qconjugate{q} a}b
  = b\qconjugate{a \qconjugate{q} a}b
  = b \qconjugate{a} q \qconjugate{a}b
  = r q s, 
  \nonumber \\
  r = b \qconjugate{a}, \qquad 
  s = \qconjugate{a}b, 
  \qquad
  \modulus{r} = \modulus{b}\modulus{a} = 1 = \modulus{s},
\end{gather} 
in two orthogonal planes\index{quaternions!rotation!four-dimensional}, exactly as studied in Section \ref{ssc:geointerp}.

The combination of three reflections at 3D hyperplanes, given by unit quaternions $a,b,c$ , leads to 
\begin{equation}
  -c\qconjugate{[-b\qconjugate{-a \qconjugate{q} a}b]}c
  = d \,\qconjugate{q} t, \quad 
  d = -c \qconjugate{b} a, \,\,\,t = a \qconjugate{b} c, \,\,\,
  \modulus{d} = \modulus{c}\modulus{b}\modulus{a} = \modulus{t} = 1.
\end{equation} 
The product of the reflection map $-\qconjugate{q}$ of \eqref{eq:3Dhypref} with $d \,\qconjugate{q} t$ leads to $-dqt$, a \emph{double rotation}\index{quaternions!rotation!double} as studied in Section \ref{sc:NewQFT}. Therefore  $d\,\qconjugate{\sandwich}\,t$ represents a \emph{rotary reflection} (rotation reflection)\index{quaternions!rotation!rotary}\index{quaternions!rotation reflection}\index{rotary reflection}\index{rotation reflection}. The three reflections $-a\qconjugate{q}a$, $-b\qconjugate{q}b$, $-c\qconjugate{q}c$ have the intersection of the three 3D hyperplanes as a resulting common pointwise invariant line, which is $d+t$, because
\begin{equation}
  d\,\qconjugate{(d+t)}\,t = d \,\qconjugate{t}t + d\,\qconjugate{d}t = d+t.
\end{equation} 
In the remaining 3D hyperplane, orthogonal to the pointwise invariant line\index{quaternions!reflection!rotary!invariant line} through the origin in direction $d+t$, the axis of the rotary reflection\index{quaternions!reflection!rotary!axis} is
\begin{equation}
  d\,\qconjugate{(d-t)}\,t = - d \,\qconjugate{t}t + d\,\qconjugate{d}t = -d+t = -(d-t).
\end{equation} 
We now also understand that a sign change of $d \rightarrow -d$
(compare three reflections at three 3D hyperplanes
$-c\qconjugate{[-b\qconjugate{(-a \qconjugate{q} a)}b]}c$
with three reflections at three lines
$+c\qconjugate{[+b\qconjugate{(+a \qconjugate{q} a)}b]}c$)
simply exchanges the roles of pointwise invariant line $d+t$ and rotary reflection axis $d-t$.

Next, we seek for an explicit description of the rotation plane of the rotary reflection $d\,\qconjugate{\sandwich}\,t$. We find that for the unit quaternions $d=e^{\alpha g}, t = e^{\beta f}$ the commutator
\be 
  [d,t] 
  = dt-td 
  = e^{\alpha g} e^{\beta f} - e^{\beta f}e^{\alpha g}
  = (gf-fg) \sin \alpha \sin \beta ,
  \label{eq:commdt}
\ee 
is a pure quaternion, because
\be 
  \qconjugate{gf-fg}= fg-gf = - (gf-fg).
\ee 
Moreover, $[d,t]$ is orthogonal to $d$ and $t$, and therefore orthogonal to the plane {spanned by} the pointwise invariant line $d+t$ and the rotary reflection axis $d-t$, because
\be
  \Scalar{[d,t]\qconjugate{d}}
  = \Scalar{dt\qconjugate{d}-td\,\qconjugate{d}}
  = 0, 
  \qquad
  \Scalar{[d,t]\qconjugate{t}} = 0 .
\ee 
We obtain a second quaternion in the plane \change{orthogonal to} $d+t$\change{, and} $d-t$\change{,} by applying the rotary reflection to $[d,t]$
\be 
  d\,\qconjugate{[d,t]}\,t
  = -d[d,t]t
  = -[d,t]\qconjugate{d}t,
\ee 
because $d$ is orthogonal to the pure quaternion $[d,t]$. We can construct an \emph{orthogonal basis of the plane of the rotary reflection}\index{quaternions!rotary reflection!rotation plane!basis} $d\,\qconjugate{\sandwich}\,t$ by computing the pair of orthogonal quaternions
\be 
  v_{1,2} 
  = [d,t]\mp d\,\qconjugate{[d,t]}\,t
  = [d,t]\pm [d,t]\qconjugate{d}t 
  = [d,t] (1\pm \qconjugate{d}t).
  \label{eq:rotrefbasis}
\ee 
For finally computing the rotation angle, we need to know the relative length of the two orthogonal quaternions $v_1, v_2$ of \eqref{eq:rotrefbasis}. For this it helps to represent the unit quaternion $\qconjugate{d}t$ as 
\be 
  \qconjugate{d}t = e^{\gamma u}, \qquad 
  \gamma \in \R, \,\,\,
  u \in \H, \,\,\, u^2=-1. 
\ee 
We then obtain for the length ratio
\begin{gather}
  r^2=\frac{|v_1|^2}{|v_2|^2}=\frac{|1+ \qconjugate{d}t|^2}{|1- \qconjugate{d}t|^2}
  = \frac{(1+ \qconjugate{d}t)(1+ \qconjugate{t}d)}{(1- \qconjugate{d}t)(1- \qconjugate{t}d)}
  = \frac{1+\qconjugate{d}t\qconjugate{t}d+\qconjugate{d}t+\qconjugate{t}d}{1+\qconjugate{d}t\qconjugate{t}d-\qconjugate{d}t-\qconjugate{t}d}
  \nonumber \\
  = \frac{2+2 \cos \gamma}{2-2 \cos \gamma}
  = \frac{1+ \cos \gamma}{1- \cos \gamma} .
\end{gather} 
By applying the rotary reflection $d\,\qconjugate{\sandwich}\,t$ to $v_1$ and decomposing the result with respect to the pair of orthogonal quaternions in the rotary reflection plane \eqref{eq:rotrefbasis} we can compute the rotation angle. 
Applying the rotary reflection to $v_1$ gives
\begin{gather}
  d\,\qconjugate{v_1}\,t
  = d \,\qconjugate{[d,t]- d\qconjugate{[d,t]}t}\,t
  = d \qconjugate{[d,t](1+\qconjugate{d}t)}t
  = d (1+\qconjugate{t}d)\qconjugate{[d,t]}t
  \nonumber \\
  = d (1+\qconjugate{t}d)(-[d,t])t
  = -[d,t] (\qconjugate{d}t + \qconjugate{d}t\qconjugate{d}t) .
\end{gather} 
The square of $\qconjugate{d}t$ is
\begin{gather}
  (\qconjugate{d}t)^2
  = (\cos \gamma + u \sin \gamma)^2
  = -1 + 2 \cos \gamma \,[\cos \gamma + u \sin \gamma]
  \nonumber \\
  = -1 + 2 \cos \gamma  \,\qconjugate{d}t .
\end{gather} 
We therefore get 
\begin{gather}
  d\,\qconjugate{v_1}\,t
  = -[d,t] (\qconjugate{d}t -1 + 2 \cos \gamma  \qconjugate{d}t)
  = [d,t] (1 - (1+ 2 \cos \gamma)  \qconjugate{d}t)
  \nonumber \\
  = a_1 v_1 + a_2 r v_2 ,
\end{gather} 
and need to solve
\be 
  1 - (1+ 2 \cos \gamma)  \qconjugate{d}t
  = a_1 (1+ \qconjugate{d}t) + a_2 r (1- \qconjugate{d}t),
\ee 
which leads to
\be 
  d\,\qconjugate{v_1}\,t
  = - \cos \gamma v_1 + \sin \gamma r v_2 
  = \cos (\pi -\gamma) v_1 + \sin (\pi -\gamma) r v_2 . 
\ee 
The rotation angle of the rotary reflection\index{quaternions!rotary reflection!rotation angle} $d\,\qconjugate{\sandwich}\,t$ in its rotation plane $v_1, v_2$ is therefore
\be 
  \Gamma = \pi -\gamma ,
  \qquad \gamma = \arccos \Scalar{\qconjugate{d}t} .
  \label{eq:Gamma}
\ee 
In terms of $d=e^{\alpha g}, t = e^{\beta f}$ we get 
\be 
  \qconjugate{d}t
  = \cos \alpha \cos \beta 
  -g \sin \alpha \cos \beta
  +f \cos \alpha \sin \beta
  -gf \sin \alpha \sin \beta .
\ee 
And with the angle $\omega$ between $g$ and $f$
\begin{gather}
  gf = \frac{1}{2}(gf+fg) +\frac{1}{2}(gf-fg)
     = \Scalar{gf} + \frac{1}{2}[g,f]
     \nonumber \\
     = -\cos \omega - \sin \omega \frac{[g,f]}{|[g,f]|},
\end{gather} 
we finally obtain for $\gamma$ the scalar part \Scalar{\qconjugate{d}t} as
\begin{gather}
  \Scalar{\qconjugate{d}t} 
  = \cos \gamma 
  = \cos \alpha \cos \beta + \cos \omega \sin \alpha \sin \beta
  \nonumber \\
  = \cos \alpha \cos \beta - \Scalar{gf} \sin \alpha \sin \beta.
  \label{eq:Sbardt}
\end{gather} 

In the special case of $g=\pm f$, $\Scalar{gf}=\mp 1$, \ie for $\omega = 0,\pi$, we get from \eqref{eq:Sbardt} that
\begin{gather}
  \Scalar{\qconjugate{d}t} 
  = \cos \alpha \cos \beta \pm \sin \alpha \sin \beta
  = \cos \alpha \cos \beta + \sin (\pm\alpha) \sin \beta
  \nonumber \\
  = \cos(\pm\alpha - \beta), 
\end{gather}
and thus using \eqref{eq:Gamma} the rotation angle\index{quaternions!rotary reflection!rotation angle} would become
\be 
  \Gamma = \pi - (\pm \alpha - \beta) = \pi \mp \alpha + \beta .
  \label{eq:Gammag=pmf}
\ee 
Yet \eqref{eq:Gamma} was derived assuming $[d,t]\neq 0$. But direct inspection shows that \eqref{eq:Gammag=pmf} is indeed correct: For $g=\pm f$ the plane $d+t, d-t$ is identical to the $1,f$ plane. The rotation plane is thus a plane of pure quaternions orthogonal to the $1,f$ plane. The quaternion conjugation in $q \mapsto d \, \qconjugate{q} \, t$ leads to a rotation by $\pi$ and the left and right factors lead to further rotations by $\mp\alpha$ and $\beta$, respectively. Thus \eqref{eq:Gammag=pmf} is verified as a special case of \eqref{eq:Gamma} for $g=\pm f$.

By substituting in Lemma \ref{lm:expqexp} $(\alpha,\beta) \rightarrow (-\beta, -\alpha)$, and by taking the quaternion conjugate we obtain the following Lemma.

\begin{lem}\label{lm:expcqexp}
Let $q_{\pm}=\frac{1}{2}(q \pm fqg)$ be the OPS of Definition \ref{df:genOPS}. For left and right exponential factors we have the identity\index{quaternions!exponential factors!identity}
\be 
  e^{\alpha g} \,\qconjugate{q_{\pm}} e^{\beta f}
  = \qconjugate{q_{\pm}} e^{(\beta\mp\alpha) f}
  = e^{(\alpha\mp\beta) g} \,\qconjugate{q_{\pm}}. 
  \label{eq:lemexpcqexp}
\ee 
\end{lem}

\subsection{New Steerable QFTs with Quaternion Conjugation and Two Pure Unit Quaternions \ensuremath{f,g}}

We therefore consider now the following new variant of the 
(double sided form of the) QFT \cite{10.1007/s00006-007-0037-8,EH:DirUP_QFT}
in $\H$ (replacing both $\i$ \change{with} $g$ and $\j$ \change{with} $f$,
and using quaternion conjugation).
It is essentially the quaternion conjugate of the new QFT of Definition \ref{df:QFTfg},
but because of its distinct local transformation geometry \change{it} deserves separate treatment. 

\begin{defn}\textbf{(QFT with respect to \ensuremath{f,g}, including quaternion conjugation)\index{quaternions!QFT!quaternion conjugation}}\label{df:QFTgf}\\
Let $f,g \in \H$, $f^2=g^2=-1$, be any two pure unit quaternions.
The quaternion Fourier transform with respect to $f,g$, involving quaternion conjugation, is 
\begin{equation}
  \label{eq:QFTgf}
  \mathcal{F}_c^{g,f}\{ h \}(\bomega)
  = \qconjugate{\mathcal{F}^{f,g}\{ h \}(-\bomega)}
  = \int_{\R^2} e^{-g x_1\omega_1} \qconjugate{h(\bvect{x})} \,e^{-f x_2\omega_2} d^2\bvect{x},
\end{equation}where $h\in L^1(\R^2,\H)$, $d^2\bvect{x} = dx_1dx_2$ and $\bvect{x}, \bomega \in \R^2 $.
\end{defn}

Linearity of the integral in \eqref{eq:QFTgf} of Definition \ref{df:QFTgf} leads to the following corrollary to Theorem \ref{th:fpmtrafo}.

\begin{cor}[QFT $\mathcal{F}_c^{g,f}$ of $h_{\pm}$]
\label{cr:gfpmtrafo}
The QFT $\mathcal{F}_c^{g,f}$ \eqref{eq:QFTgf} of the $h_{\pm} = \frac{1}{2}(h\pm fhg)$ OPS split parts\index{quaternions!QFT!quaternion conjugation!split parts}, with respect to any two unit quaternions $f,g$, of a quaternion module function 
$h \in L^1(\R^2,\H)$ have the quasi-complex forms\index{quaternions!QFT!quaternion conjugation!quasi-complex}
\begin{gather}
  \label{eq:crgfpmtrafo}
  \mathcal{F}^{g,f}_{c}\{h_{\pm}\} (\bomega)
  = \qconjugate{\mathcal{F}^{f,g}\{ h_{\pm} \}(-\bomega)}
  \nonumber \\
  \stackrel{}{=} \int_{\R^2}
    \qconjugate{h_{\pm}}e^{-f (x_2\omega_2 \mp x_1\omega_1)}d^2x
  \stackrel{}{=} \int_{\R^2}
    e^{-g (x_1\omega_1 \mp x_2\omega_2)}\qconjugate{h_{\pm}}d^2x \,\, .
\end{gather}
\end{cor}

Note, that the pure unit quaternions $f,g$ in Definition \ref{df:QFTgf} and Corrollary \ref{cr:gfpmtrafo} do not need to be orthogonal, and that the cases $f=\pm g$ are fully included. Corollary \ref{cr:gfpmtrafo} leads to discretized and fast versions of the QFT with quaternion conjugation of Definition \ref{df:QFTgf}. 

It is important to note that the roles (sides) of $f,g$ appear exchanged in   \eqref{eq:QFTgf} of Definition \ref{df:QFTgf} and in Corollary \ref{cr:gfpmtrafo}, although the same OPS of Definition \ref{df:genOPS} is applied to the signal $h$ as 
in Sections \ref{sc:OPSfg} and \ref{sc:NewQFT}. This role change is due to the presence of quaternion conjugation in Definition \ref{df:QFTgf}. Note that it is possible to first apply \eqref{eq:QFTgf} to $h$, and subsequently split the integral with the OPS$^{g,f}$ $\mathcal{F}_{c,\pm}=\frac{1}{2}(\mathcal{F}_c\pm g\mathcal{F}_cf)$, where the particular order of $g$ from the left and $f$ from the right is due to the application of conjugation in \eqref{eq:crgfpmtrafo} to $h_{\pm}$ \emph{after} $h$ is split with \eqref{eq:genOPS} into $h_+$ and $h_-$.

\subsection{Local Geometric Interpretation of the QFT with Quaternion Conjugation \label{ssc:GIntQFTc}}

Regarding the \emph{local geometric interpretation}\index{quaternions!QFT!quaternion conjugation!geometric interpretation!local} of the QFT with quaternion conjugation of Definition \ref{df:QFTgf} we need to distinguish the following cases, depending on $[d,t]$ and on whether the left and right phase factors 
\be 
  d=e^{-g x_1\omega_1}, \qquad t=e^{-f x_2\omega_2} ,
  \label{eq:dt}
\ee 
attain scalar values $\pm 1$ or not. 

Let us first assume that $[d,t]\neq 0$, which by \eqref{eq:commdt} is equivalent to $g\neq \pm f$, and $\sin(x_1\omega_1)\neq 0$, and $\sin(x_2\omega_2)\neq 0$. Then we have the generic case of a local rotary reflection with pointwise invariant line\index{quaternions!QFT!quaternion conjugation!invariant line!local} of direction 
\be
  d+t = e^{-g x_1\omega_1} + e^{-f x_2\omega_2}, 
\ee
rotation axis\index{quaternions!QFT!quaternion conjugation!rotation axis!local} in direction
\be
  d-t = e^{-g x_1\omega_1} - e^{-f x_2\omega_2}, 
\ee
rotation plane with basis \eqref{eq:rotrefbasis}, and by \eqref{eq:Gamma} and \eqref{eq:Sbardt} the general rotation angle\index{quaternions!QFT!quaternion conjugation!rotation angle!local} 
\begin{gather}
  \Gamma 
  = \pi - \arccos \Scalar{\qconjugate{d}t},
  \nonumber \\ 
  \Scalar{\qconjugate{d}t} 
  = \cos(x_1\omega_1) \cos (x_2\omega_2) 
    - \Scalar{gf} \sin (x_1\omega_1) \sin (x_2\omega_2).
\end{gather}

Whenever $g = \pm f$, or when $\sin(x_1\omega_1)=0$ ($x_1\omega_1=0,\pi[\text{mod}\, 2\pi]$, \ie $d=\pm 1$), we get for the pointwise invariant line in direction $d+t$ the simpler unit quaternion direction expression $e^{-\frac{1}{2}(\pm x_1\omega_1+x_2\omega_2)f}$, because we can apply
\be 
  e^{\alpha f}+e^{\beta f} 
  = e^{\frac{1}{2}(\alpha+\beta)f}
    (e^{\frac{1}{2}(\alpha-\beta)f}+e^{\frac{1}{2}(\beta-\alpha)f})
  = e^{\frac{1}{2}(\alpha+\beta)f} 2 \cos\frac{\alpha-\beta}{2} ,
  \label{eq:expaf+gf}
\ee 
and similarly for the rotation axis $d-t$ we obtain the direction expression
$e^{-\frac{1}{2}(\pm x_1\omega_1+x_2\omega_2+\pi)f}$, 
whereas the rotation angle is by \eqref{eq:Gammag=pmf} simply
\be
  \Gamma = \pi \pm x_1\omega_1 - x_2\omega_2. 
\ee 

For $\sin(x_2\omega_2)=0$ ($x_2\omega_2=0,\pi[\text{mod}\, 2\pi]$, \ie $t=\pm 1$), the pointwise invariant line in direction $d+t$ simplifies by \eqref{eq:expaf+gf} to $e^{-\frac{1}{2}(x_1\omega_1+x_2\omega_2)g}$, and the rotation axis with direction $d-t$ simplifies to 
$e^{-\frac{1}{2}(x_1\omega_1+x_2\omega_2+\pi)g}$, whereas the angle of rotation is by \eqref{eq:Gammag=pmf} simply 
\be
  \Gamma = \pi + x_1\omega_1 - x_2\omega_2. 
\ee

\subsection{Phase Angle QFT with Respect to \ensuremath{f,g}, Including Quaternion Conjugation}

Even when quaternion conjugation is applied to the signal $h$ we can propose a further new version of the QFT${}_c^{g,f}$, with a straightforward two \emph{phase angle} interpretation. The following definition to some degree ignores the resulting local rotary reflection effect of combining quaternion conjugation and left and right phase factors of Section \ref{ssc:GIntQFTc}, but depending on the application context, it may nevertheless be of interest in its own right.

\begin{defn}\textbf{(Phase angle QFT\index{quaternions!QFT!quaternion conjugation!phase angle transformation} with respect to  \ensuremath{f,g}, including quaternion conjugation)}\label{df:QFT2anglec}\\
Let $f,g \in \H$, $f^2=g^2=-1$, be any two pure unit quaternions. 
The phase angle quaternion Fourier transform with respect to $f,g$, involving quaternion conjugation, is 
\begin{align}
  \mathcal{F}_{cD}^{g,f}\{ h \}(\bomega)
  &= \qconjugate{\mathcal{F}_{D}^{f,g}\{ h \}(-\omega_1,\omega_2)}
  \nonumber \\
  &= \int_{\R^2} e^{-g \frac{1}{2}(x_1\omega_1+x_2\omega_2)} \qconjugate{h(\bvect{x})}
    \,e^{-f \frac{1}{2}(x_1\omega_1-x_2\omega_2)} d^2\bvect{x}.
  \label{eq:QFT2anglec}
\end{align}
where again $h\in L^1(\R^2,\H)$, $d^2\bvect{x} = dx_1dx_2$ and $\bvect{x}, \bomega \in \R^2 $. 
\end{defn}

Based on Lemma \ref{lm:expcqexp}, one possible \emph{geometric interpretation}\index{quaternions!QFT!quaternion conjugation!phase angle transformation!interpretation} of the integrand of \eqref{eq:QFT2anglec}) is a local phase rotation\index{quaternions!QFT!quaternion conjugation!phase rotation!local} of $\qconjugate{h_+}$ by angle 
$-(x_1\omega_1-x_2\omega_2)/2+ (x_1\omega_1+x_2\omega_2)/2 = + x_2\omega_2$
in the $\qconjugate{q_+}$ plane, 
and a second local phase rotation 
of $\qconjugate{h_-}$ by angle 
$-(x_1\omega_1-x_2\omega_2)/2- (x_1\omega_1+x_2\omega_2)/2 = - x_1\omega_1$ in the $\qconjugate{q_-}$ plane. This is expressed in the following corollary to Theorem \ref{th:2anglefpmtrafo}.


\begin{cor}[Phase angle QFT of $h_{\pm}$, involving quaternion conjugation]
\label{cr:fpmtrafoc}
The phase angle QFT with quaternion conjugation\index{quaternions!QFT!quaternion conjugation!phase angle transformation!split parts} of Definition \ref{df:QFT2anglec} applied to the $h_{\pm}$ OPS split parts, with respect to any two pure unit quaternions $f,g$, of a quaternion module function $h \in L^1(\R^2,\H)$ leads to the quasi-complex\index{quaternions!QFT!quaternion conjugation!phase angle transformation!quasi-complex} expressions
\begin{align}
  \mathcal{F}^{g,f}_{cD}\{h_+\}(\bomega) 
  &= \qconjugate{\mathcal{F}^{g,f}_{D}\{h_+\}(-\omega_1,\omega_2)}
  \nonumber \\
  &= \int_{\R^2}
    \qconjugate{h_{+}}e^{+f x_2\omega_2 }d^2x
  \stackrel{}{=} \int_{\R^2}
    e^{-g x_2\omega_2}\qconjugate{h_{+}}d^2x \,\, ,
\end{align}
\begin{align}
  \mathcal{F}^{g,f}_{cD}\{h_-\}(\bomega)
  &= \qconjugate{\mathcal{F}^{g,f}_{D}\{h_-\}(-\omega_1,\omega_2)}
  \nonumber \\
  &= \int_{\R^2}
    \qconjugate{h_{-}}e^{-f x_1\omega_1 }d^2x
  \stackrel{}{=} \int_{\R^2}
    e^{-g x_1\omega_1}\qconjugate{h_{-}}d^2x \,\, .
\end{align}
\end{cor} 

Note that based on Theorem \ref{th:abdetfg} the two phase rotation planes (analysis planes) are again freely steerable\index{quaternions!QFT!quaternion conjugation!phase angle transformation!steerable}. Corollary \ref{cr:fpmtrafoc} leads to discretized\index{quaternions!QFT!quaternion conjugation!phase angle transformation!discrete} and fast\index{quaternions!QFT!quaternion conjugation!phase angle transformation!fast} versions of the phase angle QFT with quaternion conjugation of Definition \ref{df:QFT2anglec}.

\section{Conclusion}

The involution maps $\i\sandwich\j$ and $f\sandwich g$ have led
us to explore a range of similar quaternionic maps\index{quaternions!quaternion maps} $q \mapsto a q b$ and $q \mapsto a \qconjugate{q} b$, where $a,b$ are taken to be unit quaternions. Geometric interpretations of these maps as reflections, rotations, and rotary reflections in 4D can mostly be found in \cite{Coxeter1946}. We have further developed these geometric interpretations to gain a \emph{complete local transformation geometric understanding}\index{quaternions!QFT!geometric understanding!local} of the integrands of the proposed new quaternion Fourier transformations (QFTs) applied to general quaternionic signals $h \in L^1(\R^2,\H)$. This new geometric understanding is also valid for the special cases of the hitherto well-known left-sided,
right-sided, and left- and right-sided (two-sided) QFTs of \cite{Ell:1993,10.1007/s00006-007-0037-8,10.1049/el:19961331,10.1109/TIP.2006.884955} and numerous other references. 

Our newly gained geometric understanding itself motivated us to propose \emph{new types} of QFTs with specific geometric properties. The investigation of these new types of QFTs\index{quaternions!QFT!new types} with the generalized form of the orthogonal 2D planes split of Definition \ref{df:genOPS} lead to important QFT split theorems\index{quaternions!QFT!split theorems}, which allow the use of discrete and (complex) Fourier transform software for efficient discretized and fast numerical implementations. 

Finally, we are convinced that our geometric interpretation of old and new QFTs paves the way for new applications, \eg, regarding \emph{steerable filter design}\index{quaternions!QFT!filter design!steerable}\index{filter design!steerable} for specific tasks in image, colour image and signal processing, \etc. 

\section*{Acknowledgements} 

E.~H. wishes to thank God for the joy of doing this research, his family, and S.~J.~S. for his great cooperation and hospitality. 

\bibliographystyle{abbrv}
\bibliography{bibliodata}

\printindex
\end{document}